\documentclass[a4paper,DIV=12,11pt,parskip=half]{scrartcl}

\usepackage[USenglish]{babel}
\usepackage[utf8x]{inputenc}
\usepackage[pdfencoding=auto,
unicode=true,
plainpages=false,
pdfpagelabels=true]{hyperref}

\usepackage{amssymb,amsmath}

\usepackage{mathtools}
\mathtoolsset{centercolon}

\usepackage[calc]{adjustbox}

\usepackage{wrapfig}
\usepackage[pangram]{blindtext}
\newlength{\strutheight}

\usepackage{graphicx}
\graphicspath{{graphics/}}

\usepackage{multirow}
\usepackage{booktabs}

\usepackage{siunitx}

\usepackage{caption}
\usepackage{subcaption}
\usepackage{amsthm}

\usepackage{enumitem}

\usepackage{cleveref}
\crefname{section}{Sec.}{Sec.}
\Crefname{equation}{Eq.}{Eqs.}
\Crefname{figure}{Fig.}{Figs.}
\Crefname{tabular}{Tab.}{Tabs.}

\usepackage{xcolor}
\usepackage[section]{placeins}

\renewcommand{\vec}{\boldsymbol}
\newcommand{\mat}{\boldsymbol}

\renewcommand{\d}{\; \mathrm d}
\renewcommand{\P}{\mathbb{P}}
\newcommand{\R}{\mathbb{R}}
\newcommand{\N}{\mathbb{N}}
\newcommand{\E}{\mathcal{E}}

\newcommand{\velocity}{\mat v}

\theoremstyle{definition}

\newtheorem{problem}{Problem}
\newtheorem{remark}{Remark}
\newtheorem{assumption}{Assumption}

\theoremstyle{remark}

\numberwithin{equation}{section}
\numberwithin{figure}{section}
\numberwithin{table}{section}
\numberwithin{problem}{section}
\numberwithin{remark}{section}
\numberwithin{assumption}{section}

\begin{document}

\title{CutFEM and ghost stabilization techniques for higher order space-time discretizations of the Navier--Stokes equations}

\author{Mathias Anselmann$^{\star,1}$, Markus Bause$^{\star}$\\
	{\small ${}^\star$ Helmut Schmidt University, Faculty of 
		Mechanical Engineering, Holstenhofweg 85,}\\ 
	{\small 22043 Hamburg, Germany}\\
}

\date{}
\maketitle

\footnotetext[1]{Corresponding author: anselmann@hsu-hh.de}

\begin{abstract}
\textbf{Abstract.} We propose and analyze computationally a new fictitious domain method, based on higher order space-time finite element discretizations, for the simulation of the nonstationary, incompressible Navier--Stokes equations on evolving domains. The physical domain is embedded into a fixed computational mesh such that arbitrary intersections of the moving domain's boundaries with the background mesh occur. The potential of such cut finite element techniques for higher order space-time finite element methods has rarely been studied in the literature so far and deserves further elucidation. The key ingredients of the approach are the weak formulation of Dirichlet boundary conditions by Nitsche's method, the flexible and efficient integration over all types of intersections of cells by moving boundaries and the spatial extension of the discrete physical quantities to the entire computational background mesh including fictitious (ghost) subdomains of fluid flow. Thereby, an expensive remeshing and adaptation of the sparse matrix data structure are avoided and the computations are accelerated. To prevent spurious oscillations caused by irregular intersections of mesh cells, a penalization, defining also implicitly the extension to ghost domains, is added. These techniques are embedded in an arbitrary order, discontinuous Galerkin discretization of the time variable and an inf-sup stable discretization of the spatial variables. The parallel implementation of the matrix assembly is described. The optimal order convergence properties of the algorithm are illustrated in a numerical experiment for an evolving domain. The well-known 2d benchmark of flow around a cylinder as well as flow around moving obstacles with arising cut cells and fictitious domains are considered further. 
\end{abstract}

\section{Introduction} 

Cut finite element methods (CutFEM) are subject to active current research. They are suitable when it comes to solving partial differential equations on evolving domains with moving boundaries.
They have been analyzed and studied numerically for parabolic problems, using low order finite elements and BDF time stepping schemes; cf. e.g. \cite{lehrenfeldEulerianFiniteElement2019}. This analysis was recently extended to the time-dependent Stokes problem using, for the spatial discretization, lowest order Taylor-Hood elements \cite{vonwahlUnfittedEulerianFinite2020} or equal-order elements along with stabilization \cite{burmanEulerianTimesteppingSchemes2020}, combined with BDF time-stepping schemes. For instance, we refer to \cite[Thm.\ 5.26]{vonwahlUnfittedEulerianFinite2020} for a rigorous error estimate.
Also, the simulation of incompressible flow on complex domains \cite{winterNitscheCutFinite2018,burmanFictitiousDomainMethods2014,schottNewFaceorientedStabilized2014,courtFictitiousDomainApproach2014} or even the simulation of coupled problems of multi-physics, for instance of fluid--structure interaction, with moving interfaces \cite{schottMonolithicApproachFluid2019,courtFictitiousDomainFinite2015,burmanNitschebasedFormulationFluidstructure2020} becomes feasible with CutFEM based techniques. The appreciable advantage of the CutFEM approach is, that all the computations are done on a time-independent, fixed background mesh. However, mesh cells then are intersected in an arbitrary manner by moving boundaries of the domain or inner interfaces, for instance, if fluid-structure interaction is studied. In contrast to the well-known Arbitrary--Lagrangian--Eulerian (ALE) method (cf.\ \cite{hughesLagrangianEulerianFiniteElement1981,doneaArbitraryLagrangianeulerianFinite1982,braessArbitraryLagrangianEulerian2000}), no transformation of the physical domain to a fixed computational reference domain is induced by CutFEM and the computational mesh is not necessarily fitted to the outer boundaries of the domain or, for  instance, the boundaries of enclosed rigid domains if flow around obstacles is studied. In classical ALE approaches, larger motions and deformations of the domains lead to a poor quality of the non-transformed, physical mesh (cf., e.g.\ \cite{wall1LargeDeformationFluidStructure2006,nazemStressIntegrationMesh2006,farinattiaymoneMeshMotionTechniques2004}), since the computational mesh has to track and resolve moving boundaries or interfaces. This mesh deformation impacts the transformation of the model equations to the reference domain and, thereby, the stability of the overall ALE approach. To preserve the mesh quality and stability of the transformation, remeshing becomes inevitable, combined with the necessity to project discrete solutions from one mesh to another, which is computationally expensive and can cause an accumulation of projection errors (cf.~\cite{brennerPrioriErrorAnalysis2014}). Fixed-mesh ALE approaches (cf.\ \cite{baigesFixedmeshALEApproach2010}) try to overcome this issue by combining a fixed background mesh with the ALE concept. In each time step, the discrete functions are projected onto a fixed background mesh, so that additional degrees of freedom, used to preserve the mesh quality by introducing supplementary grid nodes, are naturally defined by the ALE approach.

The CutFEM approach overcomes the ALE transformation of the model equations to a reference domain. In CutFEM, the model is discretized on a fixed background mesh in its native formulation. The CutFEM methodology is based on three key ingredients: Firstly, boundary conditions of Dirichlet type are imposed in a weak variational form by using Nitsche's method, cf.\ \cite{nitscheUeberVariationsprinzipZur1971,burmanFictitiousDomainFinite2012,beckerMeshAdaptationDirichlet2002,boiveauPenaltyfreeNitscheMethod2016,massingStabilizedNitscheFictitious2014,anselmannHigherOrderGalerkin2020}.
Secondly, a flexible and efficient approach to integrate over all types of mesh cells resulting from intersections of computational grid cells by moving boundaries, referred to as cut cells, is needed. The efficient integration over such cut cells requires the computation of volume integrals over portions of the underlying finite element cells.
This integration should not disturb the convergence order of the underlying space-time discretization scheme.
One way to do this is to divide the cut cells into sub-elements and to apply a standard quadrature formula on the resulting triangulation; cf.\ \cite{xieFINITEELEMENTMETHOD2011,gawlikHighorderFiniteElement2014}.
However, in this approach the mesh information has to be recomputed or rearranged for assembling the algebraic system, which counteracts the methodology and advantages of CutFEM techniques. There exist some remedies for this drawback, for instance, by using a refined sub-mesh for the integration over the cut cells, without adjusting the underlying mesh and degrees of freedom (cf.\ \cite{carraroImplementationEXtendedFinite2016,freiLocModFELocallyModified2021}), which however is still a challenging task, in particular in the case of three space dimensions. An alternative is given by using the divergence theorem and transforming the volume integrals to surface integrals (line integrals in the two-dimensional case). The latter ones can then be computed by a moment fitting method (cf.~\cite{massingEfficientImplementationFinite2013}) or suitable quadrature formulas (cf.~\cite{sudhakarAccurateRobustEasytoimplement2014}). Our approach proposes a different technique, that deviates from the previous ones. We compute iterated integrals of multivariate calculus over all portions of cells intersected by the domain's boundaries by applying (iterated) one-dimensional quadrature formulas. This results in a flexible and efficient algorithm. 

The third and last ingredient of CutFEM techniques for evolving domains is the extension of discrete functions to fictitious or ghost subdomains or the whole computational mesh, respectively, along with the ghost penalization of the discrete variational equations; cf.~\cite{burmanEulerianTimesteppingSchemes2020,vonwahlUnfittedEulerianFinite2020,burmanContinuousInteriorPenalty2007,burmanFictitiousDomainMethods2014,lehrenfeldEulerianFiniteElement2019}. The stabilization of the variational problem aims at reducing unphysical and spurious oscillations that are due to unavoidable irregular (small) cuts of mesh cells by moving boundaries. A strong analogy to the discretization of convection-dominated problems (cf.~\cite{roosRobustNumericalMethods2008}) can be observed. Further, the stabilization also aims at improving the conditioning of the resulting algebraic system; cf.~\cite{burmanGhostPenalty2010,massingStabilizedNitscheFictitious2014,lehrenfeldEulerianFiniteElement2019}. Proposed stabilization techniques are based on suitable extensions of the discrete functions to some neighborhood of the flow domain; cf.~\cite{vonwahlUnfittedEulerianFinite2020,burmanUnfittedNitscheMethod2014}. We refer to this extension of the flow domain as a fictitious or ghost subdomain, since the physical unknowns are not defined in the domain's extension by some mathematical problem. Usually, the prolongation of the discrete functions is restricted to some thin layer adjacent to the domain's boundary; cf.\ \cite{lehrenfeldEulerianFiniteElement2019,vonwahlUnfittedEulerianFinite2020,hansboCutFiniteElement2016,burmanCutFiniteElement2017}. However, in our approach we extend all quantities to the entire computational mesh. More precisely, for our prototype model problem of flow around a moving rigid obstacle in a pipe (cf.\ Fig.~\ref{fig:problem_overview}) this means that the velocity and pressure variable are extended to the whole rigid obstacle. The appreciable advantage of this augmented extension is that the underlying degrees of freedom and mesh information are fixed over the whole simulation time. This frees us from an expensive remeshing, a redistribution of degrees of freedom and a rebuilding or updating of the sparse system matrix data structure. For all time steps of the considered time interval, a non-condensed system matrix for all degrees of freedom of the whole computational grid, as the union of the flow domain and fictitious or ghost subdomains, is computed. No adaption of the matrix's data structure with respect to active degrees of freedom of the flow domain and non-active or fictitious degrees of freedom of the ghost subdomains, varying in time due to the evolving domain, is required. As concerns the stabilization itself, our approach differs from the popular derivative jump ghost penalization technique used in, e.g.,  \cite{frachonCutFiniteElement2019,hansboCutFiniteElement2016,burmanCutFEMDiscretizingGeometry2015,burmanFictitiousDomainMethods2014,schottMonolithicApproachFluid2019}. In this approach jumps in the (higher order) derivatives are penalized over facets which requires the computationally expensive evaluation of all higher order derivatives up to the maximum polynomial degree. We adopt a computationally cheaper and more direct concept of ghost penalization suggested in \cite{preussHigherOrderUnfitted2018,lehrenfeldEulerianFiniteElement2019} for linear convection-diffusion-reaction equations and generalized in \cite{vonwahlUnfittedEulerianFinite2020} to the Stokes system. In this approach, no derivatives are required and the penalization is not depending on the polynomial degree of the spatial finite element spaces, just on some scaling factors in terms of negative powers of the spatial mesh size and the viscosity. This direct stabilization is carefully adapted to the mathematical setting of the Navier--Stokes system with two unknown variables and their non-equal order, inf-sup stable approximation in space. The stability of the proposed choice of the scaling factors in terms of the spatial mesh size is demonstrated.     

Finally, we mention that an arbitrary order, discontinuous Galerkin method is applied for the discretization of the time variable. In \cite{frachonCutFiniteElement2019,zahediSpaceTimeCutFinite2017},
a piecewise linear and discontinuous in time discretization, that is combined with CutFEM techniques, is developed. The integration over the cut cells follows similar ideas as used here, but differs by subdiving the flow region of the cut cells into triangles (in two space dimensions) and integrating over these triangles whereas iterated integration over the (partially) curved subcell is applied here. Our motivation for using a discontinuous approximation in time is to overcome some difficulty in the pressure approximation of equal-order and continuous in time discretizations of Navier--Stokes solutions. Precisely, the latter lacks from the availability of an (continuous and discrete) initial value for the pressure variable that is not provided explicitly by the Navier--Stokes system, but needed for the unique definition of the full pressure trajectory. We note that discrete pressure approximations of a continuous in time approximation, for instance in the Gauss quadrature nodes of a Gauss quadrature formula in time, become accessible without any initial pressure values, which however is not sufficient to compute the complete pressure trajectory. The discontinuous discretization in time can be applied in a natural way to time-dependent domains and combined with the extension of discrete functions to fictitious domains and the integration over cut cells. By using a discontinuous test basis we reduce the space-time formulation to a time-marching scheme. Appreciable advantage of this approach is that algebraic systems of reasonable dimension are obtained. However, for higher order time discretizations the block structure of the resulting algebraic system is still complex and requires highly efficient algebraic solvers. Geometric multigrid methods have proved to be very efficient, when it comes to solving the Navier--Stokes equations \cite{hussainEfficientStableFinite2013,johnNumericalPerformanceSmoothers2000}. For non-evolving domains, a suitable preconditioning technique for the Navier--Stokes system based on a geometric multigrid method is analyzed computationally in \cite{anselmannGeometricMultigridMethod2021}.
Its extension to evolving domains and CutFEM techniques is still an ongoing research and beyond the scope of this work.    

In our numerical studies we consider a sequence of four test problems of increasing complexity to analyze the stability and performance properties of the proposed approach. Firstly, its convergence properties are analyzed for a time-dependent domain. Secondly, the well-know DFG benchmark \cite{schaferBenchmarkComputationsLaminar1996} of flow around a cylinder in two space dimensions is considered. Even though the domain is time-independent, it's worth to analyze the accuracy of the proposed approach for this well-understood setting and to compare it with discretizations based on body-fitted techniques. Thereby, our integration over cut cells, the extension of discrete functions to fictitious (ghost) subdomains of the computational mesh and the ghost penalization are evaluated. In the third and fourth experiment, time-independent domains are studied again. Flow around moving rigid bodies is investigated and illustrated for two different settings.

This work is organized as follows. In Sec.~\ref{Sec:Problem}, the prototype model problem and the notation are introduced. In Sec.~\ref{Sec:STFEM} our space-time finite element approach for simulating the Navier--Stokes system on evolving domains is presented. In Sec.~\ref{Sec:Implement}, practical aspects of the parallel implementation of the algorithms in the \textit{deal.II} library  \cite{arndtDealIILibrary2020} and the linear algebra package \textit{Trilinos} \cite{thetrilinosprojectteamTrilinosProjectWebsite2020} are addressed.      
In Sec.~\ref{sec:numerical_results}, the numerical results for a sequence of test  problems are presented. We end with a summary and an outlook in Sec.~\ref{Sec:SumOut}. 

\section{Mathematical problem and notation}
\label{Sec:Problem}

\subsection{Mathematical problem}

\begin{figure}[h!tb]
    \centering
    \includegraphics[width=0.6\linewidth]{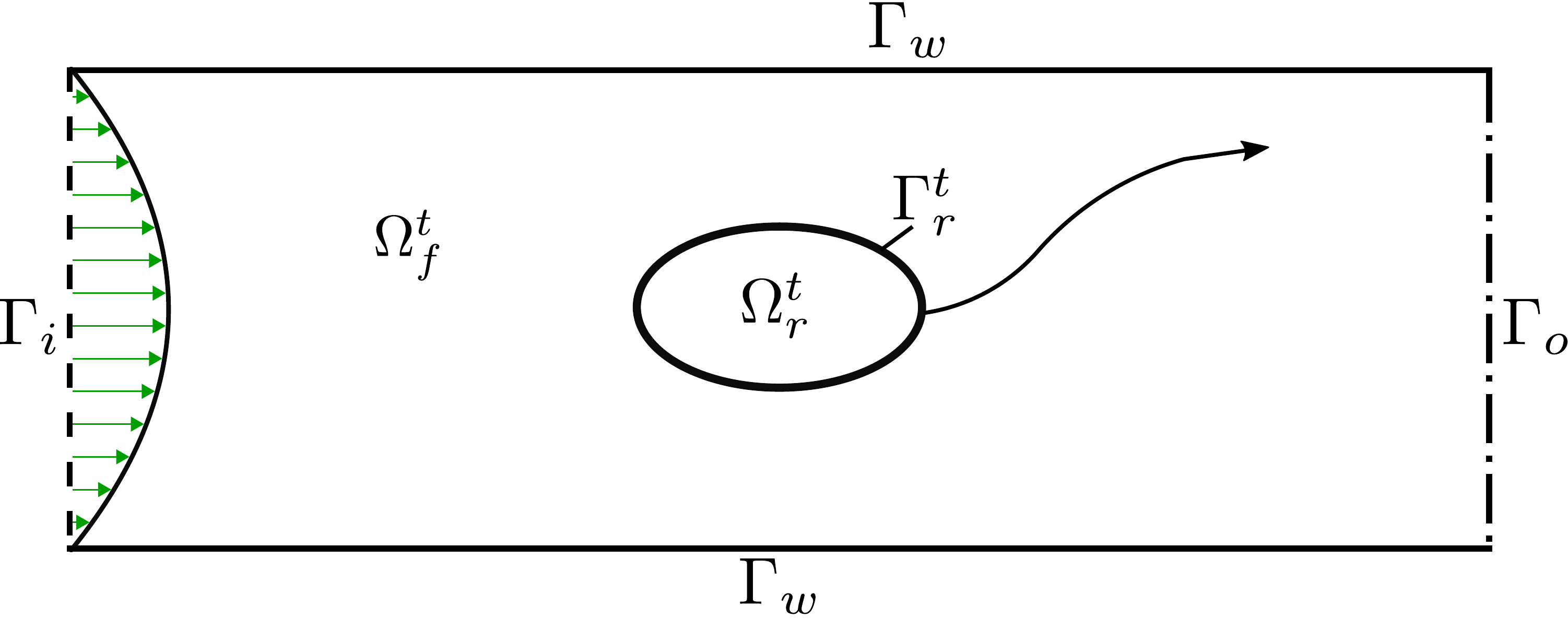}
    \caption{Problem setting and corresponding notation with pipe's boundary $\Gamma_p := \Gamma_i\cup \Gamma_w\cup \Gamma_o$.}%
    \label{fig:problem_overview}
\end{figure}

Without loss of generality, in this work we study a prototype model problem that is sketched in Fig.~\ref{fig:problem_overview}. Restricting ourselves to considering incompressible viscous flow around a moving, undeformable, rigid body in a pipe is only done in order to simplify the notation, reduce technical ballast and focus on the essential features of the numerical techniques. The problem is of high interest in practice and can be regarded as a test problem for more sophisticated applications, for instance, as a building block for fluid-structure interaction. By $\Omega=(0,L)\times (0,W)$, with some $L>0$ and $W>0$, we denote the rectangular pipe of length $L$ and width $W$. Let $T>0$ be the final time. About the rigid body $\Omega_r(t)$ we make the following assumption. 
\begin{assumption}
\label{assump:RB}
The smooth motion of the bounded rigid body $\Omega_r^t:=\Omega_r(t)$, with positive measure $\text{meas}_2(\Omega_r^t) >0$ and sufficiently smooth boundary $\Gamma_r^t :=\partial \Omega_r(t)$, for $t\in [0,T]$, is supposed to be given. Contact between the pipe boundary $\Gamma_p := \Gamma_i\cup \Gamma_w\cup \Gamma_o$ and the rigid body is assumed not to occur such that $\overline{\Omega_r^t}\subset \Omega$ and, with some suitable constant $\delta >0$,
\[
\text{dist}(\Gamma_r^t,\Gamma_p)\geq \delta > 0 \,, \quad \text{for } t\in[0,T]\,.
\]
\end{assumption}


By $\Omega_f^t:=\Omega_f(t)= \Omega\backslash \overline{\Omega_r^t}$, with $\Omega_f^t\subset \Omega$, we denote the open domain filled with fluid. For $\Omega_f^t$, with $t\in [0,T]$, and $I:=(0,T]$ we consider solving the incompressible Navier--Stokes system
\begin{subequations}
\label{eq:navier_stokes}
\begin{alignat}{4}
\label{eq:navier_stokes_0}
    \partial_t \vec{v} + (\vec{v} \cdot \nabla) \vec{v} - \nu \Delta \vec{v} + \nabla p	&= \vec{f}
    & \hspace*{2ex} & \text{in } \Omega_f^t \times I\,,
	\\
\label{eq:navier_stokes_1}
    \nabla \cdot \vec{v} &= 0
    & & \text{in } \Omega_f^t  \times I\,,
	\\
\label{eq:navier_stokes_2}
    \vec{v} & = \vec{0}
    & & \text{on } \Gamma_w \times I\,, \quad
	\\
\label{eq:navier_stokes_3}
    \vec{v} &= \vec{g}_i
    & & \text{on } \Gamma_i  \times I\,, \quad
	\\
\label{eq:navier_stokes_4}
    \vec{v} &= \vec{g}_r
    & & \text{on } \Gamma_r^t  \times I\,, \quad
	\\
\label{eq:navier_stokes_5}
    \nu \nabla \vec{v} \cdot \vec{n} - \vec{n}p &= \vec 0
    & & \text{on } \Gamma_{o}  \times I\,, \quad
	\\
\label{eq:navier_stokes_6}
    \vec{v}(0) &= \vec{v}_0 
    & & \text{in } \Omega_f^t \,.
\end{alignat}
\end{subequations}
In \eqref{eq:navier_stokes}, the velocity field $\vec v$ and the pressure $p$ are the unknown variables. In \eqref{eq:navier_stokes_0}, the parameter $\nu>0$ denotes the fluid's viscosity and the right-hand side function $\vec{f}$ is a given external force. In \eqref{eq:navier_stokes_5}, the field $\vec n$ is the outer unit normal vector. In \eqref{eq:navier_stokes_6}, the function $\vec{v}_0$ denotes the prescribed initial velocity. On $\Gamma_i$, an inflow flow profile is prescribed by \eqref{eq:navier_stokes_3}. The boundary parts $\Gamma_w$ are fixed walls with no-slip boundary \eqref{eq:navier_stokes_2}, and $\Gamma_o$ represents an outflow boundary that is modeled by the do-nothing boundary condition \eqref{eq:navier_stokes_5}; cf.~\cite{heywoodArtificialBoundariesFlux1996}. The boundary condition \eqref{eq:navier_stokes_4} ensures that the fluid follows the prescribed rigid body's motion $\vec{g}_r$ on $\Gamma_r^t$. Here, $\vec{g}_r$ is the prescribed material velocity of the particles from $\Gamma_r^t$. To simplify the notation, we put 
\begin{equation*}
	\vec{g} =
	\begin{cases} 
		\vec{0} & \text{on } \Gamma_w  \times I\,, \\
		\vec{g}_i & \text{on } \Gamma_i  \times I\,, \\
		\vec{g}_r & \text{on } \Gamma_r^t \times I
	\end{cases}
\end{equation*}
and, with $\Gamma_D^t:=\Gamma_i\cup \Gamma_w \cup \Gamma_r^t$, we rewrite the conditions \eqref{eq:navier_stokes_2} to \eqref{eq:navier_stokes_4} by 
\begin{equation*}
\label{eq:navier_stokes_7}
\vec v = \vec g \; \text{ on }  \Gamma_D^t\times I \,. 	 
\end{equation*}

We assume, that a unique (weak) solution of the system \eqref{eq:navier_stokes} exists. This implies that sufficient smoothness conditions about the evolving flow domain are satisfied. For the existence and uniqueness of solutions to the Stokes and Navier--Stokes system on time-dependent domains we refer to, e.g., \cite{burmanEulerianTimesteppingSchemes2020,salviNavierStokesEquationsNoncylindrical1988,bockNavierStokesEquationsNoncylindrical1977}. In particular, we refer to \cite[p.~156, Thm.~1, p.~164, Cor.]{salviNavierStokesEquationsNoncylindrical1988} for the existence of local in time solutions to \eqref{eq:navier_stokes} or solutions for sufficiently small data. Even though solutions to \eqref{eq:navier_stokes} might lack (due to the polygonal pipe boundary) higher order regularity, that is needed to expect an optimal order convergence behavior of higher order discretization methods, members of such families are applied here. In the numerical approximation of partial differential equations it is widely accepted that higher order methods achieve accurate results on computationally feasible grids, even if the solution of the mathematical problem is not sufficiently smooth on the whole space-time domain. This also applies to the lack of regularity of Navier--Stokes solutions for $t\rightarrow 0$ under realistic assumptions about the data of the problem (cf.~\cite{heywoodFiniteElementApproximation1982,bauseOptimalConvergenceRates2005,sonnerSecondOrderPressure2020}) which is not considered here.  

\subsection{Function spaces and forms}
\label{sec:Problem_Notation}

Here, we introduce the functions spaces that are used in this work to present our space-time CutFEM approach for the Navier--Stokes system \eqref{eq:navier_stokes}. For the analysis of solutions to \eqref{eq:navier_stokes} or the discrete scheme a more sophisticated framework is needed. The proof of well-posedness and error analyses are more involved than in the case of fixed domains. In particular, a diffeomorphism mapping the evolving domain to a reference domain and the transformation of the model equations to an equivalent system on the reference domain by the diffeomorphism are utilized strongly. For further details we refer to, e.g.,  \cite{burmanEulerianTimesteppingSchemes2020} for the Stokes system and to \cite{alphonseAbstractFrameworkParabolic2015} for an abstract framework for parabolic partial differential equations on evolving spaces. 

By $L^2(S)$ we denote the function space of square integrable functions on a domain $S$ while $H^1(S)$ is the usual Sobolev space of functions in $L^2(S)$ which have first order weak derivatives in $L^2(S)$. Further, $\langle \cdot, \cdot \rangle_S$ is the standard inner product of $L^2(S)$. We define the subspace of $L^2(S)$ with mean zero $L^2_0(S):= \{v\in L^2(S) \mid \int_S v \d x = 0\}$ and the subspace of $H^1(S)$ of functions with zero boundary values (in the sense of traces) on the Dirichlet portion $\Gamma\subset \partial S$ of the boundary $\partial S$ of $S$ as $H^1_{0,\Gamma}(S)$. Finally, by $H^{1/2}(\Gamma)$ we denote the space of all traces on $\Gamma\subset \partial S$ of functions in $H^1(S)$. For vector-valued functions we write those spaces bold. 

For the weak problem formulation we introduce the semi-linear form $A_S: (\vec H^1(S) \times L^2_0(S)) \times (\vec H^1_{0,\Gamma}(S) \times L^2_0(S)) \rightarrow \R$ by 
\begin{equation}
	\label{Eq:Slf_a}
	A_S((\vec v,p),(\vec \psi,\xi)) \coloneqq
	\langle (\vec v \cdot \vec \nabla) \vec v, \vec \psi \rangle_{S}
	+ \nu \langle \nabla \vec v , \nabla \vec \psi  \rangle_{S}
	-\langle p, \nabla \cdot \vec \psi \rangle_{S}
	+ \langle \vec \nabla \cdot \vec v, \xi \rangle_{S} 
\end{equation}
for $(\vec v, p) \in \vec H^1(S) \times L^2_0(S)$ and $(\vec \psi,\xi) \in \vec H^1_{0,\Gamma}(S) \times L^2_0(S)$. For given $\vec f \in L^2(S)$ we introduce the linear form $L: \vec H^1(S) \rightarrow \R$ by
\begin{equation}
\label{eq:DelL}
L_S(\vec \psi; \vec f) := \langle \vec f, \vec \psi \rangle_{S} 
\end{equation}
for $\vec \psi \in \vec H^1_{0,\Gamma}(S)$.

Further, we define the spaces  
\begin{equation}
\label{Eq:DefSolSpa1}
\begin{aligned}
V_{I} &\coloneqq \Big\{\vec v \in L^2\Big(I;\vec H^1\Big(\Omega_f^t\Big)\Big) \; \Big| \; \partial_t \vec v \in L^2\Big(I;\vec L^2\Big(\Omega_f^t\Big)\Big)\Big\}\,, \\
V_{0,I} &\coloneqq \Big\{\vec v \in L^2\Big(I;\vec H^1_{0,\Gamma_D^t}\Big(\Omega_f^t\Big)\Big)\Big\}
\end{aligned}
\end{equation}
and  
\begin{equation}
\label{Eq:DefSolSpa2}
L^2_{0,I} : = L^2\Big(I;L^2_0\Big(\Omega_f^t\Big)\Big)\,, \qquad L^\infty_I := L^\infty\Big(I;\vec H^1\Big(\Omega_f^t\Big)\Big)\,.
\end{equation}
In the first of the definitions in \eqref{Eq:DefSolSpa1}, the weak partial derivative with respect to the time variable $t$ is defined as an element of $\vec L^2(\Omega_{f,I})$ for the space-time flow domain $\Omega_{f,I}:= \bigcup_{0\leq t\leq T} \Omega_f^t \times \{t\}$. The definitions \eqref{Eq:DefSolSpa1} and  \eqref{Eq:DefSolSpa2} are sufficient here. For an abstract framework to treat evolving spaces we refer to \cite{alphonseAbstractFrameworkParabolic2015}. Finally, we let 
\[
\vec V_{\text{div}}^t:= \Big\{\vec v \in \vec H^1\Big(\Omega_f^t\Big) \; \Big| \; \nabla \cdot \vec v = 0\; \text{in } \Omega_f^t\Big\}\,.
\]

\subsection{Space discretization}
\label{Subsec:SpaceDisc}

\begin{figure}[!htb]
    \subcaptionbox{Computational domain $\Omega$ and background mesh $\mathcal{T}_h$. 
    \label{fig:computational_domain}}
	[0.3\columnwidth]
    {\includegraphics[height=0.2\textwidth,keepaspectratio]
    {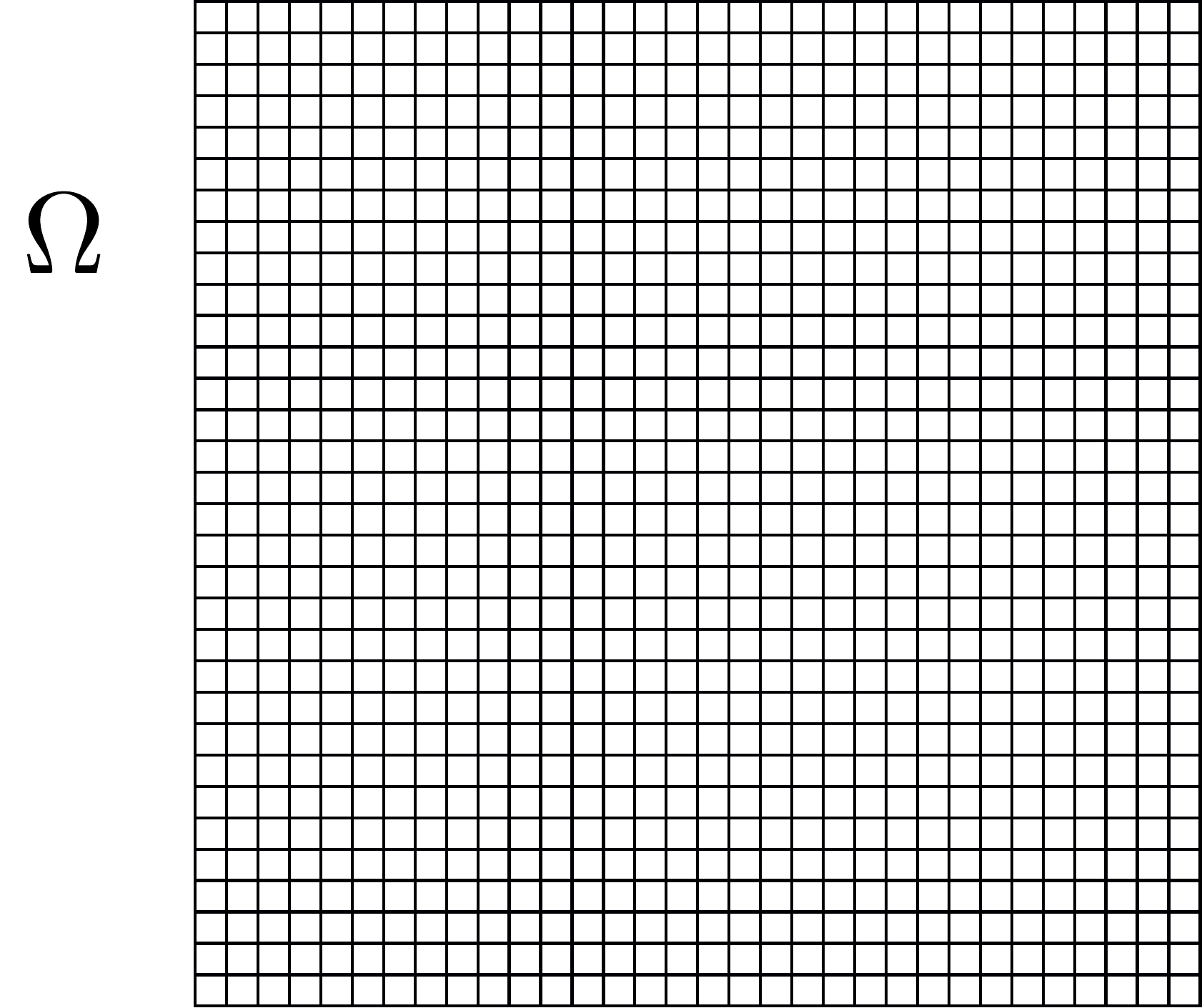}}
    \qquad
    \subcaptionbox{Fluid domain $\Omega_f^t$, rigid domain $\Omega_r^t$ and unfitted background mesh. 
    \label{fig:fluid_rigid_domain}}
	[0.28\columnwidth]
    {\includegraphics[height=0.2\textwidth,keepaspectratio]
    {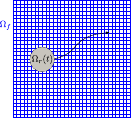}}
    \qquad
    \subcaptionbox{Submesh of fluid cells $\mathcal T_{h,f}^{t}$, rigid cells $\mathcal T_{h,r}^{t}$ and cut cells $\mathcal T_{h,c}^{t}$.
    \label{fig:cell_naming_scheme}}
    [0.28\columnwidth]
    {\includegraphics[height=0.2\textwidth,keepaspectratio]
    {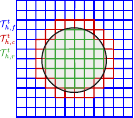}}
	\caption{CutFEM mesh topology with computational background mesh and submeshes of fluid, cut and rigid cells.}
	\label{fig:problem_cells}
\end{figure}

Let $\mathcal{T}_h=\{K\}$ be a family of shape-regular, structured decompositions of the pipe $\Omega$ (cf.\ Fig.~\ref{fig:problem_overview}) into (open) quadrilaterals $K$ with maximum cell size $h$. Precisely, $\mathcal{T}_h$ is the computational background mesh. It is not fitted to the boundary of the moving rigid body $\Omega_r^t$; cf.~Fig.~\ref{fig:problem_cells}.
We make the following assumptions.
\begin{assumption}
\label{assump:Triang}
The family $\mathcal{T}_h=\{K\}$ satisfies:
\begin{enumerate}
	\itemsep0ex
	\item The background mesh is independent of the time $t$.
	\item The mesh is cartesian, such that the cells are aligned along the coordinate lines. \label{it:cartesian}
	\item Each face of a quadrilateral is cut at most once by the rigid body \label{it:one_cut}
\end{enumerate}
\end{assumption}
The opportunity to use structured meshes facilitates their generation which is an appreciable advantage of CutFEM. The third of the assumptions is made in order to simplify the implementation when it comes to the integration over cut cells. For time-independent domains, this condition is always fulfilled, if the mesh size is chosen sufficiently small with respect to the size of the rigid body.
For evolving domains, this condition can be violated, if $\Omega_r^t$ has a curved boundary. In our implementation we ensure, that the condition of \cref{it:one_cut} is always fulfilled, cf.\ Sec.\ \ref{Sec:Implement} for details.

The time-dependent set of mesh cells $K \in \mathcal T_h$ that are subsets of the fluid domain $\Omega_f^t$ is denoted by $\mathcal T_{h,f}^{t}$, the mesh cells $K \in \mathcal T_h$ that are subsets of the rigid domain $\Omega_r^t$ is denoted by $\mathcal T_{h,r}^{t}$. The set of cut cells $K\in \mathcal T_h$ such that $K\cap \Omega_f^t \neq \emptyset$ and $K\cap \Omega_r^t \neq \emptyset$ is denoted by $\mathcal T_{h,c}^{t}$; cf.~Fig.\ref{fig:cell_naming_scheme}.

For some  $r\in \N$, let $H_h=H_h^{(r)}$ be the finite element space given by 
\begin{equation}
	\label{Eq:DefHh}
	H_h^{(r)}=\left\{v_h \in C(\overline{\Omega}) \mid v_h{}|_T\in \mathbb{Q}_r(K) \; \forall \; K \in \mathcal{T}_h \right\}\,,
\end{equation}
where $\mathbb{Q}_r(K)$ is the space defined by the multilinear reference mapping of polynomials on  the reference element with maximum degree $r$ in each variable. For the numerical experiments that are presented in Sec.~\ref{sec:numerical_results} we used the Taylor--Hood family of inf-sup stable finite element pairs for the space discretization. These elements can be replaced by any other type of inf-sup stable elements  (cf.~\cite{johnFiniteElementMethods2016}), when the proposed ghost penalty stabilization of \cref{Sec:nitsche_ghost_penalty} is applied. For brevity, we restrict ourselves to the presentation of the Taylor--Hood family. For some natural number $r\geq 2$ and with \eqref{Eq:DefHh} we then put 
\begin{equation*}
	\vec V_h = H_h^{(r)}\times H_h^{(r)} \,, \quad Q_h = H_h^{(r-1)}\,.
\end{equation*}
The space of weakly divergence free functions is denoted by
\begin{equation*}
	\vec V_h^{\text{div}} = \left\{ \vec v_h \in \vec V_h \mid \langle \nabla \cdot \vec 
	v_h,q_h\rangle  = 0 \; \text{for all } q_h \in Q_h \right\} \,.
\end{equation*}
Finally, we define the spaces  
\begin{equation*}
	\label{Eq:DefSolSpa1_h}
	V_{I,h} := \{\vec v_h \in L^2(I;\vec V_h) \mid \partial_t \vec v \in L^2(I;\vec V_h)\}\,, \qquad L^2_{0,I,h} : = L^2(I;Q_h)\,.	
\end{equation*}

\subsection{Time discretization}

For the time discretization, we decompose the time interval $I=(0,T]$ into $N$ subintervals $I_n=(t_{n-1},t_n]$, $n=1,\ldots,N$, where $0=t_0<t_1< \cdots < t_{N-1} < t_N = T$ such that $I =\bigcup_{n=1}^N I_n$ and $I_n \cap I_m = \emptyset$ for $n \neq m,\,m,n = 1,\ldots,N$.
We put $\tau = \max_{n=1,\ldots, N} \tau_n$ with $\tau_n = t_n-t_{n-1}$. Further, the set $\mathcal{M}_\tau := \{I_1,\ldots, I_N\}$ of time intervals is called the time mesh. For a Banach space $B$ of functions defined on the time-independent domain $\Omega$ and any $k\in \N_0$, we let 
\begin{align}
	\label{Def:PkIB}
	\P_k(I_n;B)
	= \Big\{w_\tau : I_n \to B \mid w_\tau(t) = \mbox{$\sum\limits_{j=0}^k$}
	W_j \, t^j \; \forall t \in I_n\,, \;
	W_j \in B\; \forall j \Big\}\,.
\end{align}
For an integer $k\in \N_0$, we put 
\begin{equation}
	\label{Def:Xtauk}
	X_\tau^{k} (B)
	\coloneqq
	\left\{w_\tau \in L^2(I;B) \mid w_\tau|_{I_n} \in
	\P_{k}(I_n;B)\; \forall I_n\in \mathcal{M}_\tau\,, \; w_\tau (0) \in B \right\}\,. 
\end{equation}

\section{Space--time finite element discretization with CutFEM}
\label{Sec:STFEM}

Here we introduce our approximation of the Navier--Stokes system \eqref{eq:navier_stokes} by combining CutFEM and variational time discretization techniques. Firstly, we introduce the variational formulation of the Navier--Stokes system \eqref{eq:navier_stokes}. The next subsection is devoted to the discretization in space. Nitsche's method \cite{beckerMeshAdaptationDirichlet2002} is applied to incorporate Dirichlet boundary conditions in a weak form. Further, the ghost penalty stabilization along with the extension of the discrete in space functions to fictitious subdomains of fluid flow, i.e.\ to the domain of the moving rigid body here, are introduced. Together, the stabilization and the extension provide an implicit definition of the discrete functions on the background mesh for the time-independent pipe domain $\Omega$. The stabilization reduces spurious oscillations by irregular cuts of finite elements by the moving rigid body $\Omega_r^t$ and improves the conditioning of the resulting algebraic system; cf.~\cite{lehrenfeldEulerianFiniteElement2019}. Finally, the discretization in time by the discontinuous Galerkin method is done and the resulting fully discrete problem is presented. In this work, the space-time finite element approach is recovered as a time marching scheme by the choice of a discontinuous test basis. This is done with the perspective of getting linear systems of reasonable size, even though the systems still have a complex block structure due to the higher order time discretization; cf.~\cite{anselmannHigherOrderGalerkin2020}.

\subsection{Variational formulation of the continuous problem}

A sufficiently regular solution of the Navier--Stokes system \eqref{eq:navier_stokes} satisfies the following variational problem that provides the basis for space-time finite element discretizations.  

\begin{problem}[Continuous problem]
\label{Prob:WNS}
Let $\vec f \in L^2(\Omega_{f,I})$ and $\vec v_0 \in \vec V_{\text{div}}^0$ be given. Let $\vec{\hat{g}} \in \vec V_{I}\cap L_I^\infty$ denote a prolongation of the boundary values on $\Gamma_D^t$ to $\Omega_{f}^t$, such that $\vec{\hat{g}} = \vec g$ on $\Gamma_D^t\times I$. Find $(\vec v,p)\in (V_I\cap L_I^\infty)\times  L_{0,I}^2$, with $\vec v \in \vec{\hat{g}} + V_{0,I}$, such that $\vec v(0) = \vec v_0$ and for all $\vec{\phi} \coloneqq (\psi,\xi) \in V_{0,I}\times L_{0,I}^2$, 
\begin{equation}
\label{eq:var_problem}
\int_0^T \langle\partial_t \vec v, \vec \psi \rangle_{\Omega_f^t} +
A_{\Omega_f^t}((\vec v,p),(\vec \psi,\xi)) \d t = \int_0^T
L_{\Omega_f^t}(\vec \psi; \vec f) 
\d t\,.
\end{equation}
\end{problem}

\begin{remark} Regarding the well-posedness of Problem~\ref{Prob:WNS} we note the following. 
\begin{itemize}
\item In \cite[p.~156, Thm.~1]{salviNavierStokesEquationsNoncylindrical1988}, the existence and uniqueness of local in time weak solutions to Eqs.~\eqref{eq:navier_stokes_0}, \eqref{eq:navier_stokes_1}, \eqref{eq:navier_stokes_6} equipped with homogeneous Dirichlet boundary conditions is shown under sufficient regularity assumptions about the data and motion of the evolving domain. The technical details of the assumptions about the domain motion are skipped here. In particular, $\vec v \in L^\infty_I$ and $\partial_t \vec v \in \vec L^2(\Omega_{f,I})$ is ensured on a sufficiently small time interval. Further, in \cite[p.~164, Cor.]{salviNavierStokesEquationsNoncylindrical1988} the existence of weak solutions is proved under smallness assumptions about the data. 

\item For evolving domains, the functional analytical setting for proving the existence of weak solutions differs from the one used for time-independent domains. Stricter assumptions about the data and domain are required. Then, a velocity field with time derivative $\partial_t \vec v\in \vec L^2(\Omega_{f,I})$ instead of  the distribution $\partial_t \vec v \in L^2(I;H^{-1}(\Omega))$ for time-independent domains $\Omega$ is obtained. For this reason, the stricter assumption $\vec f \in \vec L^2(\Omega_{f,I})$ about the right-hand side function $\vec f$ is made here. Further regularity conditions about $\vec f$ are supposed in \cite{salviNavierStokesEquationsNoncylindrical1988} for the proof of existence and uniqueness. 

\item For $\vec v \in V_I\cap L^\infty_I$ and $p\in L^2_{0,I}$ the left-hand side terms in \eqref{eq:var_problem} are well-defined. The existence of the sufficiently smooth extension $\vec{\hat{g}}$ is tacitly assumed here. This requires smoothness assumptions about the evolving boundary of the fluid domain and assumptions about polygonal portions of the boundary. We skip the technical details of the assumptions since we focus on the discrete scheme here. On the discrete finite element level we avoid such extensions by employing Nitsche's method. 
\end{itemize}  
\end{remark}

\subsection{Semidiscretization in space with weak enforcement of boundary conditions by Nitsche's method and ghost penalty stabilization}
\label{Sec:nitsche_ghost_penalty}

To incorporate Dirichlet boundary conditions on unfitted meshes we use Nitsche's method and follow \cite{beckerMeshAdaptationDirichlet2002,benkNitscheMethodNavier2012}. Here, the Dirichlet boundary conditions are enforced in a weak form by adding face integrals to the variational problem.  We also refer to \cite{anselmannHigherOrderGalerkin2020} where Nitsche's method is applied to a higher order variational time discretization of the Navier--Stokes system on time-independent domains. Degrees of freedom assigned to  Dirichlet portions of the boundary are now treated as unknowns of the variational problem. They are no longer enforced by the underlying function space and their implementation into the algebraic system, as it is usually done, with subsequent condensation of the algebraic equations.

Let $t\in [0,T]$. For the treatment of Dirichlet boundary conditions by Nitsche's method (cf.~\cite{beckerMeshAdaptationDirichlet2002}) we introduce the bilinearform  $B_{\Gamma_D^t} : \vec H^{1/2}(\Gamma_D^t) \times (\vec V_h \times Q_h) \rightarrow \R$ by 
\begin{equation}
	\label{Def:B_GamD}
	\begin{aligned}
		B_{\Gamma_D^t}(\vec w,(\vec \psi_h,\xi_h)) : = &  - \langle \vec w, \nu \nabla \vec \psi_h \cdot \vec n + \xi_h  \vec n  \rangle_{\Gamma_D^t}
		\\[1ex]
		& \quad + \gamma_1 \nu \langle h^{-1} \vec w , \vec 
		\psi_h  \rangle_{\Gamma_D^t}  + \gamma_2 \langle h^{-1} \vec w \cdot \vec n, \vec \psi_h \cdot \vec n 
		\rangle_{\Gamma_D^t} 
	\end{aligned}
\end{equation}
for $\vec w \in \vec H^{1/2}(\Gamma_D^t)$ and $(\vec \psi_h,\xi_h) \in \vec V_h \times Q_h$, where $\gamma_1>0$ and $\gamma_2> 0$ are numerical (tuning) parameters for the penalization. In  \cite{schottNewFaceorientedStabilized2014,agerNitschebasedCutFinite2019,winterNitscheCutFinite2018}, it is demonstrated computationally that their choice in the interval $(10, 100)$ leads to robust results.
In our simulation presented in Sec.~\ref{sec:numerical_results}, we put $\gamma_1 = \gamma_2 = 35$. The discrete semilinear form $A_h: (\vec V_h\times Q_h)\times (\vec V_h\times Q_h)\rightarrow \R$ is then given by 
\begin{equation}
	\label{Def:Ah}
	A_h((\vec v_h,p_h),(\vec \psi_h,\xi_h)) :=  A_{\Omega_f^t}((\vec v_h,p_h),(\vec \psi_h,\xi_h)) - \langle \nu \nabla \vec 
	v_h \cdot \vec n - p_h \vec n, \vec \psi_h\rangle_{\Gamma_D^t} + B_{\Gamma_D^t}(\vec v_h,\vec 
	\phi_h)
\end{equation}
for $(\vec v_h,p_h)\in \vec V_h\times Q_h$ and $(\vec \psi_h,\xi_h) \in \vec V_h\times Q_h$. The linear form $L_h: (\vec V_h\times Q_h) \rightarrow \R$ is defined by
\begin{equation}
\label{Def:Lh}
L_h((\vec \psi_h,\xi_h);\vec f, \vec g)	:= L_{\Omega_f^t}(\vec \psi_h; \vec f) + B_{\Gamma_D^t}(\vec g,(\vec \psi_h,\xi_h)) 
\end{equation}
for $(\vec \psi_h,\xi_h) \in \vec V_h\times Q_h$. In \eqref{Def:Ah} and \eqref{Def:Lh}, the forms $A_h$ and $L_h$ defined in  \eqref{Eq:Slf_a} and \eqref{eq:DelL}, respectively, are extended naturally to test functions of $\vec V_h$ with nonhomogeneous Dirichlet conditions.

\begin{remark}
\begin{itemize}	
\item We note that the forms \eqref{Def:B_GamD} and \eqref{Def:Ah} are introduced for (spatially) discrete functions that are defined on the time-independent background mesh $\mathcal T_h$ of the pipe $\Omega$. In \eqref{Def:B_GamD} and \eqref{Def:Ah}, the integration is done over the evolving fluid domain $\Omega_f^t$ and its boundary part $\Gamma_D^t$ where Dirichlet boundary conditions are prescribed. The domain of the rigid body $\Omega_r^t$ is considered as a fictitious (ghost) flow domain. In the following, the discrete fluid velocity $\vec v_h$ and pressure $p_h$ are defined implicitly in the ghost domain $\Omega_r^t$ by the (semi-) discrete variational problem that is augmented by a ghost penalty stabilization exploiting an extension of the discrete functions to $\Omega_r^t$. 
\item We comment on the different boundary terms  in the forms \eqref{Def:B_GamD} and \eqref{Def:Ah}.
The second term on the right-hand side of \eqref{Def:Ah} reflects the natural boundary condition, making the weak imposition of the boundary conditions consistent.
The first term on the right-hand side of \eqref{Def:B_GamD} is introduced to preserve the symmetry
properties of the continuous system. The
last two terms are penalizations, that ensure the stability of the discrete system. In the inviscid
limit $\nu = 0$, the last term amounts to a "no-penetration" condition. Thus, the form \eqref{Def:B_GamD}
provides a natural weighting between boundary terms corresponding to viscous effects $(\velocity = \vec{g})$,
convective behavior $((\velocity \cdot \vec{n})^- \velocity = ( \vec{g} \cdot \vec{n})^- \vec{g})$ and inviscid behavior $(\velocity \cdot \vec{n} = \vec{g} \cdot \vec{n})$.
\end{itemize}
\end{remark}

\begin{figure}[htpb]
	\centering
	\includegraphics[width=0.25\linewidth]{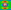}
	\caption{A moving rigid domain at two different time points $\tilde t_{n-1}$ and $\tilde t_{n}$ with computational background mesh.}%
	\label{fig:why_extension}
\end{figure}

\begin{figure}[h!tb]
	\centering
	\subcaptionbox{Irregular tiny cut cells of the background mesh due to the non body-fitted triangulation.
    \label{fig:gp_setting}}
	[1.\columnwidth]
	{
		\hspace*{-4em}
		\includegraphics[width=0.55\columnwidth,keepaspectratio]{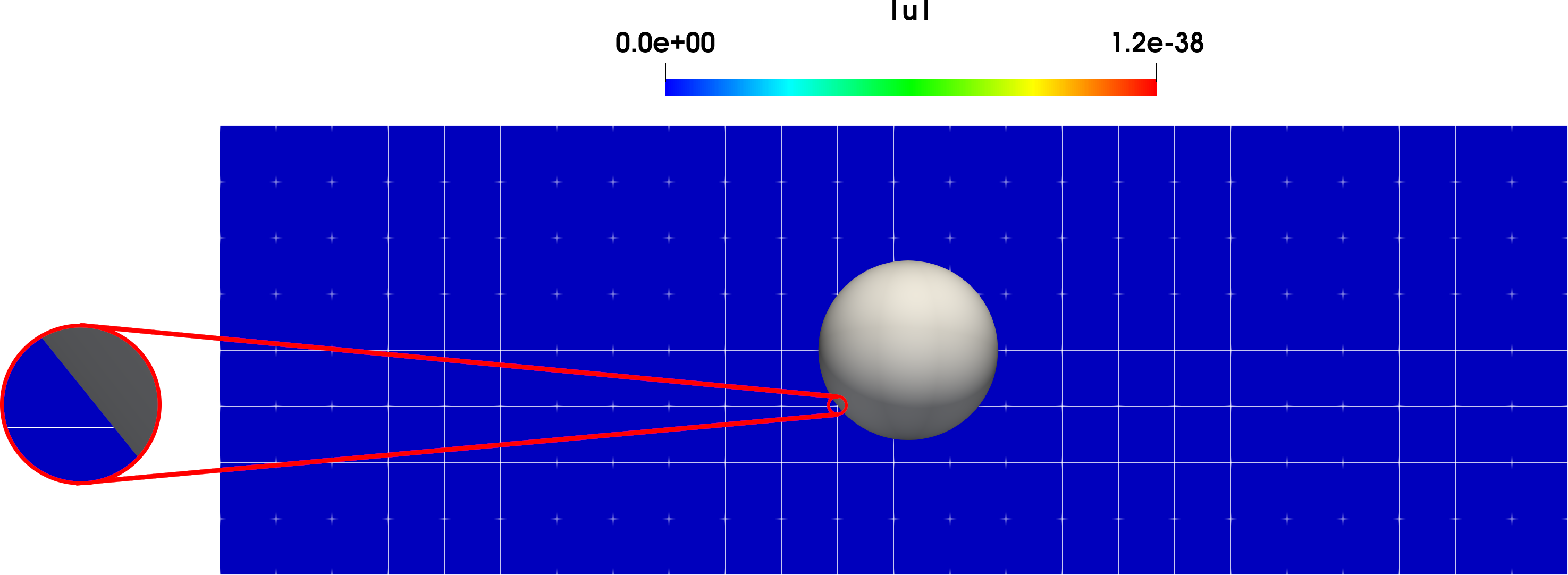}
	}
	\\[2ex]
	\subcaptionbox{Solution without ghost-penalty stabilization and extension.
		Cells $\in T_{h,r}^{t}$ are colored black and were inactive.\label{fig:gp_no}}
	[0.45\columnwidth]
	{\includegraphics[width=0.45\columnwidth,keepaspectratio]{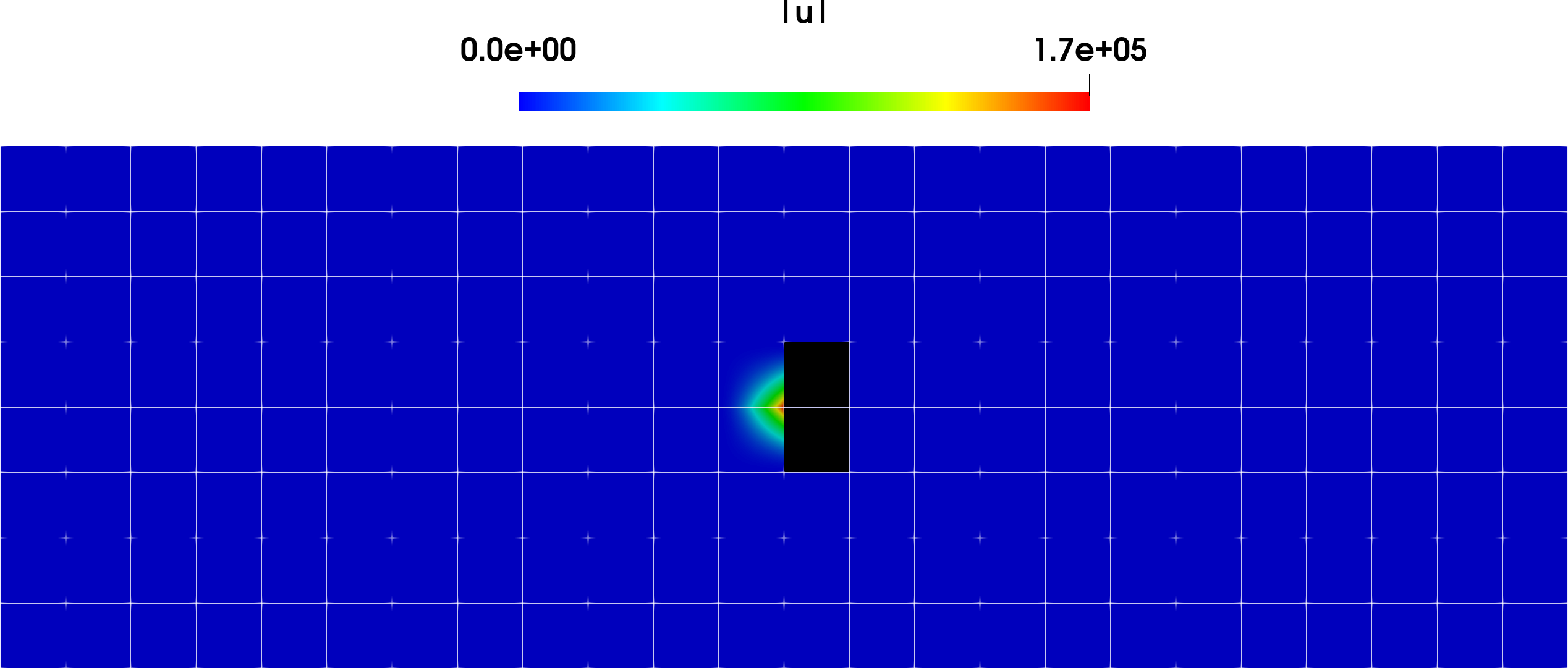}
	}
	\hspace*{3ex}
	\subcaptionbox{Solution with ghost-penalty stabilization and fictitious domain extension to rigid body $\Omega_r^t$. \label{fig:gp_yes}}
	[0.45\columnwidth]
	{\includegraphics[width=0.45\columnwidth,keepaspectratio]{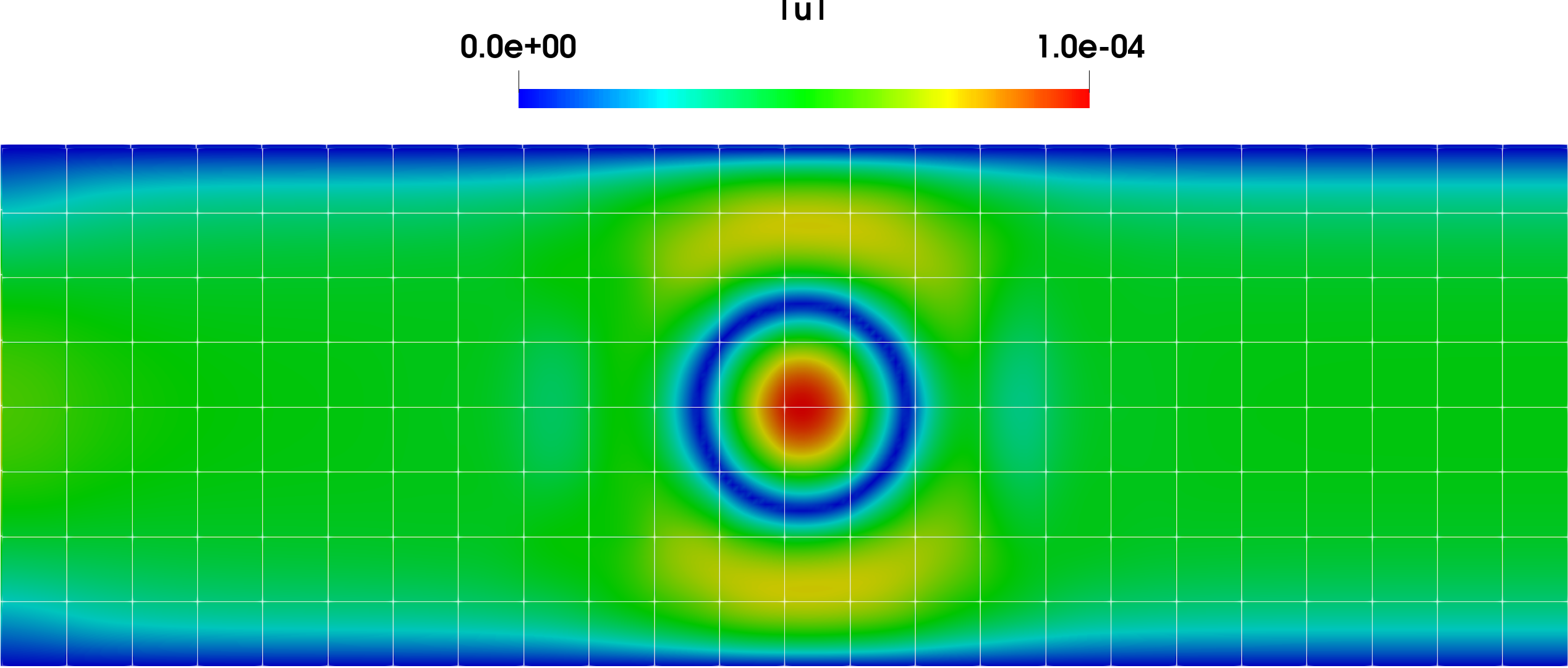}
	}
	\caption{Effect of ghost-penalty stabilization and extension.}%
	\label{fig:why_gp}
\end{figure}

To control irregular cuts of finite elements by the evolving domain and extend the fluid velocity and pressure functions from $\Omega_f^t$ to $\Omega_r^t$, and thereby to the entire, time-independent background domain $\Omega$, we introduce a discrete ghost penalty operator $S_h$. This stabilization then expands the semi-linear form \eqref{Def:Ah} by additional terms. Our stabilization modifies an approach that was firstly introduced in \cite{lehrenfeldEulerianFiniteElement2019} for convection-diffusion-reaction equations and, then, considered in \cite{vonwahlUnfittedEulerianFinite2020} for Stokes problems to the Navier--Stokes system \eqref{eq:navier_stokes}. This approach offers the appreciable advantage over further implementations of the ghost penalty method (cf.~\cite{lehrenfeldEulerianFiniteElement2019}) that no higher order spatial derivatives have to be computed. Thereby, it leads to reduced computational costs. The twofold motivation of the ghost penalty stabilization is sketched in Fig.~\ref{fig:why_extension} and Fig.~\ref{fig:why_gp}, respectively. The first aim of the ghost stabilization $S_h$, defined in \eqref{eq:extension_operator}, is to extend the solution from the physical fluid domain $\Omega_f$ to the rigid domain $\Omega_r$. This is illustrated in Fig.~\ref{fig:why_extension}. To assemble the algebraic system, i.e.\ the Jacobian matrix and right-hand-side function of the Newton linearization that is applied here, the discrete solution that is already computed for some time point $\tilde t_{n-1}$ in the corresponding fluid domain $\Omega_f^{\tilde t_{n-1}}$ needs to be evaluated in the fluid domain $\Omega_f^{\tilde t_n}$ at time $\tilde t_n$. For instance, in the case of the lowest-order discontinuous Galerkin time discretization the time points $\tilde t_{n-1}$ and $\tilde t_n$ correspond to the midpoints of the subintervals $I_{n-1}$ and $I_n$ such that $\tilde t_{n-1}=(t_{n-2}+t_{n-1})/2$ and $\tilde t_n=(t_{n-1}+t_n)/2$. In the case of higher order discontinuous Galerkin time discretizations, more time nodes or degrees of freedom in time are involved, but the extension problem applies similarly. Due to the motion of the rigid body $\Omega_r^t$, the discrete velocities and pressure values at time $\tilde t_{n-1}$ are not necessarily defined in those parts of the rigid domain $\Omega_r^{\tilde t_{n-1}}$ that belong to the fluid domain $\Omega_f^{\tilde t_{n}}$ at time $\tilde t_{n}$ such that an appropriate extension of the solution at time $\tilde t_{n-1}$ to the rigid domain at time $\tilde t_{n-1}$ or to the entire time-independent domain $\Omega$, as done in this work, becomes indispensable. For this reason, we denote the domain of the rigid body as the fictitious or ghost domain of fluid flow, since discrete velocities and pressures are not defined here by means of a physical problem, but by some artificial extension only.    

The second aim of the stabilization is to reduce spurious and unphysical oscillations. They are due to unavoidable irregular or tiny cuts of finite elements by the moving rigid domain. This is illustrated in Fig.~\ref{fig:gp_setting}. Here, we refer to the color range indicating the instability. The irregular cuts lead to steep gradients and oscillations in the computed solution profile. This is observation is related to the well-known challenges in the numerical approximation of convection-dominated transport. By steep gradients in the solution profile, the condition number of the corresponding linear algebraic system increases strongly which prevents their efficient solution by iterative methods. This problem of instability is sketched in Fig.~\ref{fig:gp_no}. For the sake of simplicity, flow around a fixed (non-moving) cylinder is simulated on a computational background mesh with arising cut cells. Here, no extension of the discrete solution to the rigid body is applied. In cut cells, the integration is done over the fluid part of the cell only. In Fig.~\ref{fig:gp_no}, an instability of the discrete solution in the cut cells is observed. By usage of ghost penalty stabilization that is proposed in the following, a smooth extension is obtained; cf.\ Fig.~\ref{fig:gp_yes}. We explicitly note, that this extension admits no reasonable physical interpretation, since it is not based on any mathematical model but represents a numerically motivated approach.

\begin{figure}[h!tb]
	\centering
	\subcaptionbox{Canonical extension $\E v_1$ of a Lagrange basis function from $K_1$ to $\omega_F=K_1\cup K_2$. \label{fig:extension_two_cells}}
	[0.38\columnwidth]
	{\includegraphics[width=0.45\columnwidth,keepaspectratio]{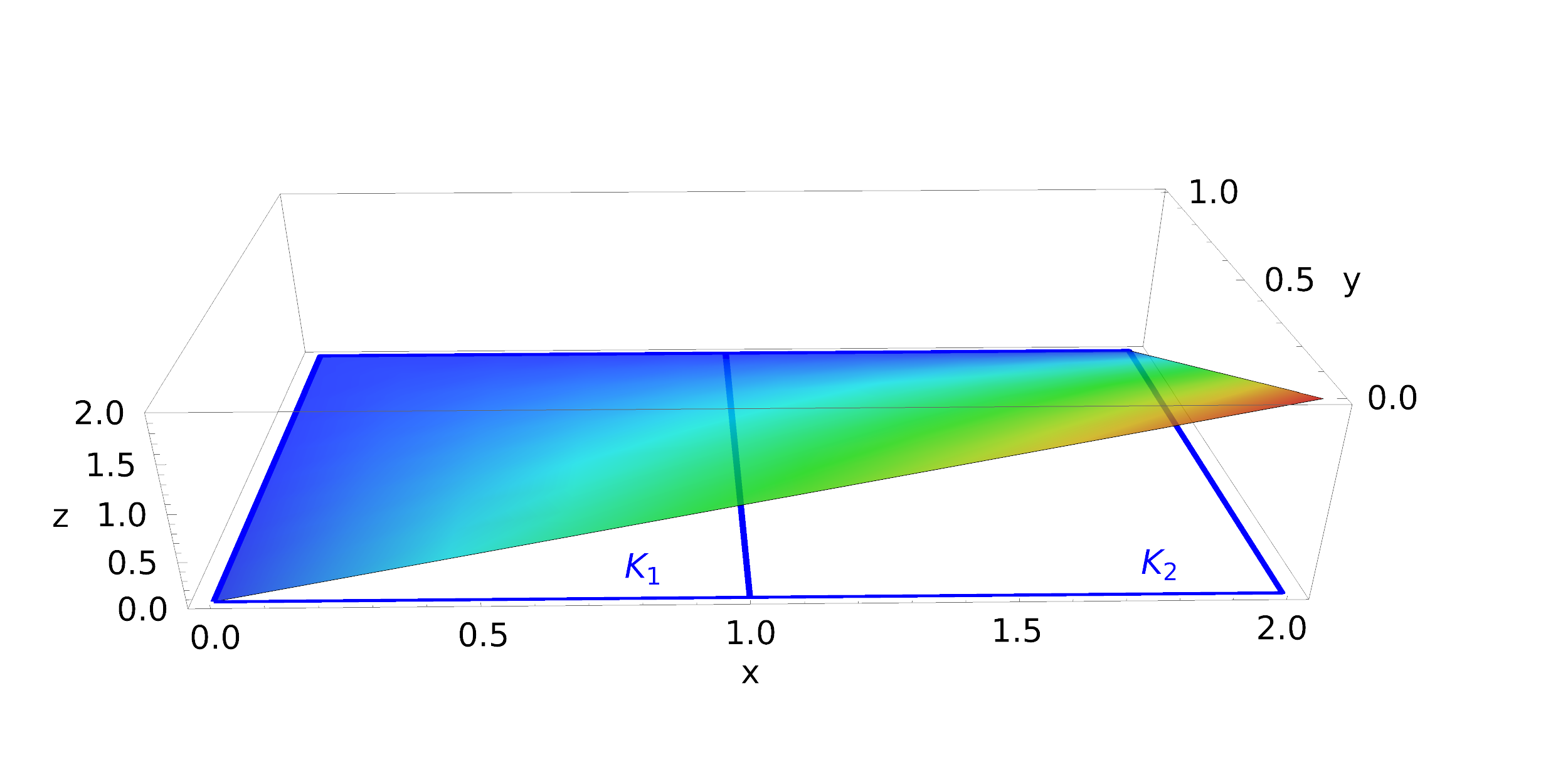}
	}
    \hspace*{7ex}
	\subcaptionbox{Patchwise extension to a fictitous domain for a cut cell $K_1$, with tiny, irregular cut.. \label{fig:extension_fictitious_domain}}
    [0.45\columnwidth]
    {\includegraphics[width=0.45\columnwidth,keepaspectratio]{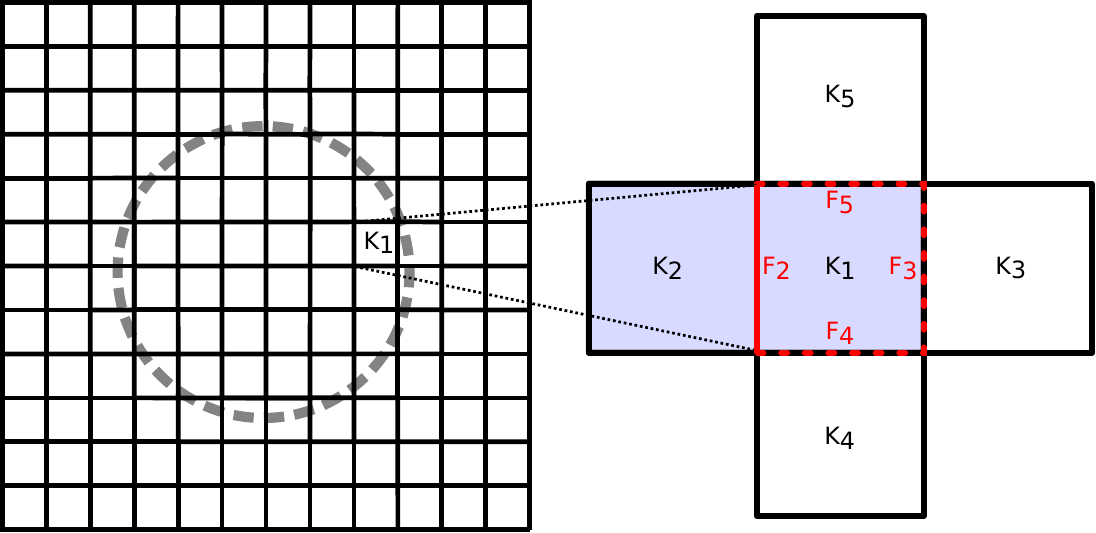}
    }	
\caption{Patchwise extension of ghost penalty stabilization.}%
\end{figure}

In the sequel, we extend the ghost penalty approach that is proposed in \cite{lehrenfeldEulerianFiniteElement2019} for linear transport equations and extended in \cite{vonwahlUnfittedEulerianFinite2020} to Stokes problems to the case of the Navier--Stokes system. However, our stabilization differs from the one that is suggested in \cite{vonwahlUnfittedEulerianFinite2020}. In \cite{vonwahlUnfittedEulerianFinite2020} the stabilization covers a boundary strip, i.e.\ the cut cells and cells of the fluid and the rigid domain adjacent to the boundary $\Gamma_r^t$. In the analysis, this is exploited to prove the error estimates. Here, we apply the combined stabilization and extension on the submeshes of the cut cells and the cells covering the rigid body (fictitious domain) $\Omega_r^t$ and, thus, to the entire domain $\Omega$ of our problem setting (cf.\ Fig.~\ref{fig:problem_overview}). This increases the efficiency of the parallel implementation (cf.\ Sec.~\ref{Sec:Implement}). Our numerical experiments (cf.\ Sec.\ref{sec:numerical_results}) indicate the admissibility of this modified ghost penalty stabilization and extension. For the problem setting of Fig.~\ref{fig:problem_overview} we denote by the time-dependent submesh $\mathcal T_{h,s}^{t}$ the set of all cut cells and all cells that are entirely in the rigid domain, 
\[
\mathcal T_{h,s}^{t} := \mathcal T_{h,r}^{t}\cup \mathcal T_{h,c}^{t} = \left\{K\in \mathcal T_h \mid K \subset \Omega_r^t \;\text{or}\; K\cap \Omega_{f}^t \neq \emptyset \neq K\cap \Omega_{r}^t \right\}\,.
\]
The set of all the faces that are common to two cells $K_1$ and $K_2$ of the submesh $\mathcal T_{h,s}^{t}$ is defined by 
\begin{equation}
	F_h^t \coloneqq \left\{
	\overline{K_1} \cap \overline{K_2} \mid 
	K_1 \in \mathcal T_{h,s}^{t}\,,\; 
	K_2 \in \mathcal T_{h,s}^{t}\,,\;
	K_1 \neq K_2\,,\; \text{meas}_1(\overline{K_1} \cap \overline{K_2})> 0
	\right\}
	\,.
\end{equation}
For $F\in F_h^t$, we denote by $\omega_F = K_1 \cup K_2$ the patch of the two adjacent elements $K_1$ and $K_2$ with the common face $F$. Then we define the bilinear form $S_{F_h^t}: (\vec V_h \times Q_h) \times (\vec V_h \times Q_h) \rightarrow \R$ by 
\begin{multline}
	\label{eq:extension_operator}
	S_{F_h^t}((\vec v_h,p_h),(\psi_h,\xi_h))
	=
	\sum_{F \in F_h^t}
	\gamma_v \left(\frac{1}{\nu} + \nu\right) \frac{1}{h^2}
	\left<(\E \vec{v}_{h|K_1} - \E \vec{v}_{h|K_2}),
	(\E \vec{\psi}_{h|K_1} - \E \vec{\psi}_{h|K_2})
	\right>_{\omega_F}
	\\
	+ \gamma_p \frac{1}{\nu}
	\left<
	(\E p_{|K_1} - \E p_{|K_2}),
	(\E \xi_{|K_1} - \E \xi_{|K_2}) 
	\right>_{\omega_F} \,,
\end{multline}
with numerical parameters $\gamma_v, \gamma_p > 0$. The stabilization \eqref{eq:extension_operator} is based on \cite{vonwahlUnfittedEulerianFinite2020}. However, in this work it is applied to the cut and rigid body cells given by $\mathcal T_{h,s}^{t}$. This approach differs from \cite{vonwahlUnfittedEulerianFinite2020} which is motivated by our aim to use a time-independent background mesh for computational efficiency. The performed computations (cf.\ Sec.~\ref{sec:numerical_results}) illustrate the robustness of the extension and stabilization. The restriction of the ghost penalty stabilization to a strip around the fluid domain is studied in a further work; cf.\ \cite{anselmannNumericalConvergenceDiscrete2021}. A superiority of the latter stabilization is not observed. In Eq.~\eqref{eq:extension_operator}, $\E$ is the canonical extension of a polynomial function of degree $r$ in each variable, that is defined on one element of the face patch $\omega_F$, to the whole patch; cf.~\cite{lehrenfeldEulerianFiniteElement2019}. Thus, we have that   
\begin{equation}
\label{def:CanExt}
	\E : \mathbb Q_r(K_i) \to \mathbb Q_r(\omega_F)\,, \quad \text{for} \; i \in \{1,2\}\,, 
\end{equation}
with $(\E u)_{|K_i} = u$ for $u\in \mathbb Q_r(K_i)$. In the case of vector-valued functions, the extension \eqref{def:CanExt} is applied component wise. To illustrate the extension \eqref{def:CanExt}, Fig.~\ref{fig:extension_two_cells} shows the extended function $\E v_1$ of a scalar valued Lagrange $\mathcal Q_1$ finite element Lagrangian basis function $v_1$ that is defined on the quadrilateral $K_1$ and has the value 1 in the grid node $(1,0)^\top$. In Fig.~\ref{fig:extension_fictitious_domain}, we illustrate the application of the stabilization operator \eqref{eq:extension_operator} for the cut cell $K_1$. This cell has four faces, that all belong to the set $F_h^t$.
For each of these faces we apply the stabilization \eqref{eq:extension_operator} along with the extension \eqref{def:CanExt}.
We start with the face $F_2$, marked by a full red line. The corresponding face patch, built from the cells $K_1$ and $K_2$, i.e.\ $\omega_F=K_1\cup  K_2$, is colored in purple in Fig.~\ref{fig:extension_fictitious_domain}.
After evaluating the extension \eqref{eq:extension_operator} on this first face patch we continue with the remaining three face patches, that are respectively built from the faces marked by the dashed red line. Proceeding in this way with all faces in $F_h^t$, we implicitly extend a discrete velocity $\vec{v}_h$ or pressure $p_h$ to $\Omega_r^t$ by adding the stabilization \eqref{eq:extension_operator} to the space discretization. The implicit definition of the extension of the discrete functions comes through the fact that the function values in the extended (ghost) domain $\Omega_r^t$ are obtained by the solution of the algebraic system and is not explicitly prescribed.  

\noindent Finally, we define the stabilized semilinear form $A_h^s: (\vec V_h\times Q_h)\times (\vec V_h\times Q_h)\rightarrow \R$ by 
\begin{equation}
	\label{Eq:AhS}
	A_h^s((\vec v_h,p_h),(\vec \psi_h,\xi_h)) := A_h((\vec v_h,p_h),(\vec \psi_h,\xi_h)) + S_{F_h^t}((\vec v_h,p_h),(\psi_h,\xi_h))
\end{equation}
for $(\vec v_h,p_h)\in \vec V_h\times Q_h$ and $(\vec \psi_h,\xi_h) \in 
\vec V_h\times Q_h$.

\begin{remark}
We note that $A_h^s$ is defined on the time-independent bulk space $\vec V_h\times Q_h$. The volume integrals of $A_h$ (cf.~\eqref{Def:Ah}) in \eqref{Eq:AhS} and of $L_h$ in \eqref{Def:Lh} are computed over the time-dependent fluid domain $\Omega_f^t$. Thus, an integration over the fluid portion of the cut cells has to be provided. This is done in Sec.~\ref{Sec:Implement}. The penalty stabilization $S_h$ (cf.~\eqref{eq:extension_operator}) is computed over time-dependent patches that are subsets either of the fluid and the rigid domain or of the rigid one only.
\end{remark}
\begin{remark}
Since we extend the solution of the physical fluid domain $\Omega_f$ to the whole computational domain $\Omega$, the choice of the time step size is in principle not restricted by a geometric CFL condition as it arises in \cite{schottMonolithicCutFinite2019}. This is confirmed by numerical experiments. Nevertheless, the choice of the time step size affects the performance and robustness of the ghost penalty stabilization. Its sensitivity depends on the geometrical and flow parameters like the diameter of the rigid domain and the fluid viscosity such that the time step size should be chosen reasonably. 
\end{remark}

We are now in a position to define a semidiscrete approximation of the system~\eqref{eq:navier_stokes}.

\begin{problem}
\label{Prob_SemiDisc}
Let $\vec f \in L^2(\Omega_{f,I})$ and $\vec v_{0,h} \in \vec V^{\text{div}}_h$ be given. Find $(\vec v_h,p_h)\in V_{I,h} \times L_{0,I,h}^2$, such that $\vec v_h(0) = \vec v_{0,h}$ and for all $(\vec \psi,\xi) \in L^2(I;\vec V_h)\times L_{0,I,h}^2$,
	\begin{equation}
		\label{eq:var_problem_space_discrete}
		\int_0^T \langle\partial_t \vec v_h, \vec \psi_h \rangle_{\Omega_f^t} +
		A_h^s((\vec v_h,p_h),(\vec \psi_h,\xi_h)) \d t = \int_0^T
		L_h(\psi_h; \vec f, \vec g) 
		\d t\,.
	\end{equation}
\end{problem}

\subsection{Fully discrete problem}
For the discretization in time we use the discontinuous Galerkin method with piecewise polynomials in time of order $k\in \N_0$. We note that continuous in time (or even continuously differentiable in time) variational discretizations of partial differential equations are known to be more efficient if the number degrees of freedom in time is measured versus the convergence rate of the time discretization if the underlying basis functions in time and quadrature formulas are chosen properly. For this observation we refer to, e.g. \cite{hussainEfficientStableFinite2013,hussainHigherOrderGalerkin2011,kocherVariationalSpaceTime2014}. The superiority of the continuous in time families is then due the fact that some degrees of freedom in time are directly obtained by continuity constraints and do not have to be computed as parts of the algebraic system; cf.~\cite{anselmannHigherOrderGalerkin2020,anselmannNumericalStudyGalerkin2020}. However, if the Stokes or Navier--Stokes system is considered, the application of continuous in time variational discretization becomes more involved. This is due to the fact that no initial value for the pressure variable is given by the mathematical model. However, such initial value is required for the unique definition of the pressure trajectory as long as equal order in time discretizations of the velocity and pressure variable are desired. Applying the discontinuous in time Galerkin method, frees us from an initial value for the pressure. Moreover, stronger stability properties of the discretization are ensured since discontinuous Galerkin methods are known to be strongly $A$-stable.

In order to keep this work self-contained, we briefly present a formal derivation of the discontinuous Galerkin discretization for the abstract evolution problem 
\begin{equation}
\label{Eq:AbstEvolProb}
\partial_t \vec v + \vec A \vec v = \vec f\,,  \qquad \vec v(0)= \vec v_0\,,
\end{equation}
as an equality in the dual space $V^\ast(t)$ of a time-dependent Hilbert space $V(t)$ of functions defined on an evolving domain $\Omega^t$ for $t\in [0,T]$; cf.~\cite{alphonseAbstractFrameworkParabolic2015}. In the formal derivation we tacitly assume that the solution of \eqref{Eq:AbstEvolProb} satisfies all the additional conditions that are required such that the arising terms are well-defined. For the derivation of the discontinuous Galerkin method we use the Reynold's transport theorem that reads as
\begin{equation}
	\label{eq:reynolds}
	\frac{\d}{\d t} \int_{\Omega^t} \vec{v} \cdot \vec{\psi} \, \d \vec{x}
	=
	\int_{\Omega^t} \partial_t \vec{v} \cdot \vec{\psi} \, \d \vec{x}
	+ \int_{\Omega^t} \vec{v} \cdot \partial_t \vec{\psi} \, \d \vec{x}
	+ \int_{\partial \Omega^t} (\vec{v} \cdot \vec{\psi}) (\vec{w} \cdot \vec{n}) \, \d s \,.
\end{equation}
Here, the vector field $\vec w(\vec x, t)$ is the material velocity of the particles from the boundary $\partial \Omega^t$, and $\vec n$ is the outer unit normal vector. Substituting \eqref{Eq:AbstEvolProb} into \eqref{eq:reynolds} yields that 
\begin{equation}
	\int_{0}^{T}
	\frac{\d}{\d t}
	\langle\vec v, \vec \psi  \rangle_{\Omega^t} 
	-
	\langle\vec v, \partial_t \vec \psi  \rangle_{\Omega^t} 
	-
	\big \langle \vec v \cdot \vec \psi,\vec w \cdot \vec n\big \rangle_{\partial \Omega^t} 
	+ \langle \vec A \vec v,\vec \psi  \rangle_{\Omega^t} 
	\d t
	= \int_{0}^{T} \langle \vec f,\vec \psi  \rangle_{\Omega^t}  
	\d t \,.
\end{equation}
With the fundamental theorem of calculus and for test functions $\psi$ with $\vec \psi(T) = \vec 0$ we get that 
\begin{equation}
\label{Eq:DervDG_1}
	- \langle\vec v(0), \vec \psi(0)  \rangle_{\Omega^0}
	-
	\int_{0}^{T}
	\langle\vec v, \partial_t \vec \psi  \rangle_{\Omega^t} 
	+
	\big \langle \vec v \cdot \vec \psi, \vec w \cdot \vec n \big \rangle_{\partial \Omega^t} 
	- \langle \vec A \vec v,\vec \psi  \rangle_{\Omega^t}  
	\d t
	= 
	\int_{0}^{T}
	\langle \vec f,\vec \psi  \rangle_{\Omega^t}  
	\d t \,.
\end{equation}
Rewriting the integrals in \eqref{Eq:DervDG_1} as a sum over the subintervals of the time mesh $\mathcal M_\tau$, we get that 
\begin{equation}
	\label{Eq:DervDG_2}
	- \langle\vec v(0), \vec \psi(0)  \rangle_{\Omega^0}
	-\sum_{n = 1}^{N}
	\int_{t_{n-1}}^{t_n}
	\langle\vec v, \partial_t \vec \psi  \rangle_{\Omega^t} 
	+
	\big \langle \vec v \cdot \vec \psi, \vec w \cdot \vec n \big \rangle_{\partial \Omega^t} 
	- \langle \vec A \vec u,\vec \psi  \rangle_{\Omega^t}  
	\d t
	= 
	\sum_{n = 1}^{N}
	\int_{t_{n-1}}^{t_n}
	\langle \vec f,\vec \psi  \rangle_{\Omega^t}  
	\d t \,.
\end{equation}
From \eqref{eq:reynolds} we derive that 
\begin{equation*}
\begin{aligned}
\label{Eq:DervDG_3}
- \int_{t_{n-1}}^{t_n} \langle\vec v, \partial_t \vec \psi  \rangle_{\Omega^t}  \d t & = - \int_{t_{n-1}}^{t_n} \frac{\d}{\d t} \langle\vec v, \vec \psi   \rangle_{\Omega^t} \d t + \int_{t_{n-1}}^{t_n} \langle \partial_t \vec v, \vec \psi  \rangle_{\Omega^t}  \d t + \int_{t_{n-1}}^{t_n} 
\big \langle \vec v \cdot \vec \psi, \vec w \cdot \vec n \big \rangle_{\partial \Omega^t} \d t\,.
\end{aligned}
\end{equation*} 
Applying the fundamental theorem of calculus to the first term on the right-hand side of this identity and substituting the resulting equation into \eqref{Eq:DervDG_2} shows that 
\begin{multline}
\label{Def:DG}
	\sum_{n = 1}^{N}
	\int_{t_{n-1}}^{t_n}
	\langle \partial_t \vec v, \vec \psi  \rangle_{\Omega^t}	
	+ \langle A\vec v,\vec \psi \rangle_{\Omega^t}
	\d t +
	\sum_{n = 1}^{N-1}
	\bigg [ \big \langle \vec v, \vec \psi  \big \rangle_{\Omega^t} \bigg ]_n
	+ \big \langle \vec v(0^+), \vec \psi(0)  \big \rangle_{\Omega^0}
	\\[1ex]
	= 
	\sum_{n = 1}^{N} \int_{t_{n-1}}^{t_n}
	\langle \vec \phi ; \vec f\rangle_{\Omega^t}
	\d t
	+
	\langle\vec v_0, \vec \psi(0) \rangle_{\Omega^0}
	 \,,
\end{multline}
with the jump operator $[\cdot]_n$ at $t_n$, defined as $[g]_n \coloneqq g(t_n^+) - g(t_n^-)$ with the one-sided limits $g(t_n^\pm) = \lim_{t\rightarrow t_n^\pm} g(t)$. For a smoothly evolving domain $\Omega^t$, that is assumed here, the one-sided limits of $\Omega^t$ at the time nodes $t_n$ coincide. We note that \eqref{Def:DG} continues to be well-defined for functions that are differentiable piecewise in time (with respect to the time mesh $\mathcal M_\tau$) only. Thus, \eqref{Def:DG} can be solved within the space \eqref{Def:Xtauk} of piecewise polynomials in time.  

Applying this concept of discontinuous Galerkin time discretization to the semidiscrete Problem~\ref{Prob_SemiDisc} yields the following fully discrete problem. 
\begin{problem}
\label{Prob_FullyDisc}
Let $\vec f \in L^2(\Omega_{f,I})$ and $\vec v_{0,h} \in \vec V^{\text{div}}_h$ be given. For $n=1,\ldots, N$, and given $\vec v_{\tau,h}{}_{|I_{n-1}} \in \mathbb P_k(I_{n-1}; \vec V_h )$ for $n>1$ and $\vec v_{\tau,h}{}_{|I_{n-1}}(t_{n-1}^-) := \vec v_{0,h}$ for $n=1$, find $(\vec v_{\tau,h},p_{\tau,h})\in \mathbb P_k(I_n; \vec V_{h}) \times \mathbb P_k(I_n;Q_h)$, such that for all $(\vec \psi_{\tau,h},\xi_{\tau,h}) \in \mathbb P_k(I_n;\vec V_h)\times \mathbb P_k(I_n;Q_h)$,
\begin{multline}
\label{Prob_FullyDisc_1}
	\int_{t_{n-1}}^{t_{n}} \langle\partial_t \vec v_{\tau,h}, \vec \psi_{\tau,h} \rangle_{\Omega_f^t} +
	A_h^s((\vec v_{\tau,h},p_{\tau,h}),(\vec \psi_{\tau,h},\xi_{\tau,h})) \d t + \langle \vec v_{\tau,h}(t_{n-1}^+), \vec \psi _{\tau,h}(t_{n-1}^+) \rangle_{\Omega^{t_{n-1}}}  \\ = \int_{t_{n-1}}^{t_{n}} 
	L_h(\psi_{\tau,h}; \vec f, \vec g) \d t + \langle \vec v_{\tau,h}(t_{n-1}^-), \vec \psi _{\tau,h}(t_{n-1}^+) \rangle_{\Omega^{t_{n-1}}}\,.
\end{multline}
\end{problem}

\begin{remark}
\label{Rem:AlgSolv}
\begin{itemize}
\item In Problem~\ref{Prob_FullyDisc} a smoothly evolving domain is assumed such that $\Omega^{t_{n-1}^+}=\Omega^{t_{n-1}^-}$ is satisfied. All volume integrals that arise in the forms of Eq.~\eqref{Prob_FullyDisc_1} are computed over the fluid domain $\Omega_f^t$.


\item To solve the algebraic counterpart of Eq.~\eqref{Prob_FullyDisc_1} we use an inexact Newton method. To enhance the range of convergence of the standard Newton method, a linesearch is applied to damp the length of a Newton step. Further, a "dogleg approach" (cf., e.g. \cite{pawlowskiInexactNewtonDogleg2008}), that belongs to the class of trust-region methods and offers the advantage that also the search direction, not just its length, can be adapted to the nonlinear solution process, was implemented. Both schemes require the computation of the Jacobian matrix of the algebraic counterpart of Eq.~\eqref{Prob_FullyDisc_1}. In the dogleg method multiple matrix-vector products with the Jacobian matrix have to be computed. Since the Jacobian matrix is stored as a sparse matrix, this matrix--vector product can be computed at low computational costs.
From the point of view of convergence, both methods yield a superlinear convergence behavior. In our numerical examples of Sec.~\ref{sec:numerical_results}, both modifications of Newton's method lead to comparable results. We did not observe any convergence problems for these nonlinear solvers.

\item To solve the linear system of the Newton iteration, we use the parallel, sparse direct solver {SuperLU\_DIST} \cite{liSuperLUDISTScalable2003}. The development of an efficient geometric multigrid preconditioner for the CutFEM approach is still an ongoing work. For our geometric multigrid preconditioning technique for space-time finite element discretizations of the Navier--Stokes problem on fitted spatial spaces for time-independent domains we refer to \cite{anselmannGeometricMultigridMethod2021}.
\end{itemize}
\end{remark}

\subsection{Algebraic in time formulation}

In this section we derive a semi-algebraic formulation of Problem~\ref{Prob_FullyDisc}. The algebraic counterpart of Eq.~\eqref{Prob_FullyDisc_1} is presented with respect to the time variable only. The restriction to the time-discretization is done due to the challenges related to the evolving domain. The presentation of the fully algebraic counterpart of Eq.~\eqref{Prob_FullyDisc_1} is skipped here for brevity. The transformation of the space discretization into an algebraic form follows the usual steps of finite element methods. For the derivation of the nonlinear algebraic system of a similar higher order space-time finite element approach to the Navier--Stokes system on time-independent domains and the application of Newton's method for the system's linearization we refer to \cite{anselmannHigherOrderGalerkin2020}. 

\begin{figure}[h!tb]
	\centering
	\includegraphics[width=0.5\linewidth]{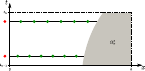}
	\caption{Space-time quadrature nodes for a one-dimensional problem and a single cut cell with a rigid domain $\Omega_r^t$.}%
	\label{fig:space_time_slab}
\end{figure}

In Problem~\ref{Prob_FullyDisc}, numerical integration in time by the $m$-point Gauss quadrature is still applied. For $g \in \{f: (0,T)\rightarrow L^2(\Omega) \mid f\in C(I_n;L^2(\Omega))\; \forall I_n \in \mathcal M_\tau\}$ this formula is given by 
\begin{equation}
\label{Def:GaussQuad}
Q_n(g) = \sum_{\mu = 1}^{m} w_\mu g(t_{n,\mu}) \approx \int_{I_n} g(t) \d t\,,
\end{equation} 
where $t_{n,\mu} \in I_n$, for $\mu = \ldots, m$, are the integration points and $w_\mu>0$ the weights. The number $m\in \N$ of quadrature  points, and thereby the order of exactness $2m-1$ the quadrature formula \eqref{Def:GaussQuad}, is chosen such that $m\geq (3k+1)/2$ is satisfied. Then, the integration in time of the discrete semilinearform $A_h((\vec v_{\tau,h},p_{\tau,h}),(\vec \psi_{\tau,h},\xi_{\tau,h}))$ on the left-hand side of \eqref{Prob_FullyDisc_1} is done exactly.
In Fig.~\ref{fig:space_time_slab} we illustrate the distribution of  the space-time integration points for a one-dimensional setting in space and a single cell. Two quadrature nodes are used for the time domain. This corresponds to the application of the dG(1) scheme in time. Six quadrature nodes are used for spatial integration, which arises in the spatial approximation by piecewise quadratic functions.

Further, we define a temporal basis $\{\chi_l\}_{l=0}^k\subset \mathbb P_k(I_n;\R)$ by the conditions 
\[
\chi_l\Big(t^{\text{GR}}_{n,\mu}\Big) = \delta_{l,\mu}\,, \quad \text{for} \;\; l,\mu =0,\ldots, k\,,
\]
for the Gauss--Radau quadrature nodes $\{t^{\text{GR}}_{n,\mu}\}_{\mu=0}^k$ of the subinterval $I_n$. Expanding the discrete solution $(\vec{v}_{\tau}{}_|{}_{I_n}, p_{\tau}{}_|{}_{I_n}) \in \mathbb P_k  (I_n;\vec {V_h}) \times \mathbb P_k (I_n;Q_h)$ on $I_n$ in terms of the basis functions yields 
\begin{align}
	\label{eq:time_discrete_ansatz}
	\vec{v}_{\tau,h|I_n}(\mat x, t)
	=
	\sum_{l=0}^{k}
		\vec {v}_{n,l} (\mat{x})\chi_{l}(t)  \quad \text{and} \quad 
	p_{\tau,h|I_n}(\mat x, t)
	=
	\sum_{l=0}^{k}
	p_{n,l} (\vec{x})\chi_{l}(t) \,,
\end{align}
for $t\in I_n$ and with coefficient functions $\vec{v}_{n,l} \in \vec V_{h}$ and $p_{n,l} \in Q_{h}$. Appreciable advantage of the usage of the Gauss--Radau quadrature points for the construction of the temporal basis is that the quantity $\vec{v}_{\tau,h|I_{n-1}}(\vec{x}, t_{n-1}^-)$ in \eqref{Prob_FullyDisc_1}, that is due to the discontinuous Galerkin scheme, is a degree of freedom of the time discretization and, thus, is available without further computational costs. 

Then, Problem~\ref{Prob_FullyDisc} can be recovered in the following form. 

\begin{problem}
	\label{prob:Algeb_time_disc}
	Let $\vec f \in L^2(\Omega_{f,I})$ and $\vec v_{0,h} \in \vec V_{h}^{\text{div}}$ be given. For $n=1,\ldots, N$, and given $\vec v_{\tau,h}{}_{|I_{n-1}} \in \mathbb P_k(I_{n-1}; \vec V_{h} )$ for $n>1$ and $\vec v_{\tau,h}{}_{|I_{n-1}}(t_{n-1}^-) := \vec v_{0,h}$ for $n=1$, find $(\vec v_{\tau,h},p_{\tau,h})\in \mathbb P_k (I_n; \vec V_{h}) \times \mathbb P_k (I_n; Q_{h})$, such that 
	for all $j \in \{0 \ldots k\}$ and $(\vec \psi_{h},\xi_{h}) \in \vec{V}_{h}\times Q_{h}$,
	\begin{multline}
		\label{eq:Algeb_timeDisc}
		\sum_{\mu=1}^{m} w_\mu \biggl[
		\langle
		\partial_t 	\vec{v}_{\tau,h|I_n}(t_{n,\mu}),
		\vec{\psi}_{h} \, \chi_j(t_{n,\mu})
		\rangle_{\Omega_f^t}
		+
		A_h^s\left(\left(
		\vec{v}_{\tau,h|I_n}(t_{n,\mu}),
		p_{\tau,h|I_n}(t_{n,\mu})
		\right),
		\left(
		\vec{\psi}_{h} \,\chi_j(t_{n,\mu}),
		\xi_{h} \,\chi_j(t_{\mu_i})
		\right) \right) \biggr]
		\\
		+
		\left\langle
		\vec{v}_{\tau,h|I_n}(t_{n-1}^+),
		\vec{\psi}_{h}\,\chi_j(t_{n-1}^+)
		\right\rangle_{\Omega^{t_{n-1}}}
		=
		\sum_{\mu=1}^{m} w_\mu \left[ \left(
		L_h(\vec{\psi}_{h}\, \chi_j(t_{n,\mu}); \vec f, \vec g)
		\right) \right]
		+
		\\
		\left\langle \vec{v}_{\tau,h|I_{n-1}}(\vec{x}, t_{n-1}^-),
		\vec{\psi}_{h}(\vec{x}) \chi_j(t_{n-1}^+)
		\right\rangle_{\Omega^{t_{n-1}}}\,.
	\end{multline}
\end{problem}

To illustrate the structure of Eq.~\eqref{eq:Algeb_timeDisc}, we review the terms in Eq.~\eqref{eq:Algeb_timeDisc} more in detail. Exemplarily, this is done for some of them. Using the expansions given in \eqref{eq:time_discrete_ansatz}, we obtain that 
\begin{align}
	\langle
	\partial_t 	\vec{v}_{\tau,h|I_n}(t_{n,\mu}),
	\vec{\psi}_{h} \, \chi_j(t_{n,\mu})
	\rangle_{\Omega_f^t}
	& =
	\sum_{l=0}^{k}
	\langle
	\vec{v}_{n,l} \, \partial_t \chi_l(t_{\mu_i}),
	\vec{\psi}_{h}\, \chi_j(t_{\mu_i})
	\rangle_{\Omega_f^t}\,,
	\\
    \langle \nabla \vec{v}_{\tau,h|I_n}(t_{n,\mu}),
	\nabla \vec{\psi}_{h}\, \chi_j(t_{\mu_i})
	\rangle_{\Omega_{f}^t}
	& =
	\sum_{l=0}^{k}
	\langle \nabla \vec{v}_{n,l} \, \chi_l(t_{n,\mu}),
	\nabla \vec{\psi}_{h} \, \chi_j(t_{n,\mu})
	\rangle_{\Omega_{f}^t}\,,
	\\
	\langle
	p_{\tau,h|I_n}(t_{n,\mu}),
	\nabla \cdot \vec{\psi}_{h}\, \chi_j(t_{n,\mu})
	\rangle_{\Omega_{f}^t}
	& =
	\sum_{l=0}^{k}
	\langle
	p_{n,l} \, \chi_{l}(t),
	\nabla \cdot \vec{\psi}_{h}\, \chi_j(t_{n,\mu})
	\rangle_{\Omega_{f}^t}\,.
\end{align}
Thus, Problem~\ref{prob:Algeb_time_disc} aims at computing the coefficient functions $\{\vec{v}_{n,l},p_{n,l}\}_{l=0}^k$, with $\{ \vec{v}_{n,l},p_{n,l}\}\in \vec V_h\times Q_h$ for $l=0,\ldots,k$. Clearly, the pair $\{\vec{v}_{n,l},p_{n,l} \} \in \vec V_h\times Q_h$ yields the fully discrete approximation in the Gauss--Radau quadrature node $t_{n,l}\in I_n$. The coefficient functions are elements of the bulk finite element spaces $\vec V_h$ and $Q_h$, respectively, and are thus defined on the computational background mesh $\mathcal T_h$ of $\Omega$. In the rigid body domain they are defined implicitly by means of the stabilization and extension operator $S_{F_h^t}$ introduced in Eq.~\eqref{eq:extension_operator}.

\section{Implementational aspects}
\label{Sec:Implement}

In this section we address key aspects of the implementation of the presented CutFEM higher order space-time approach with arbitrary polynomial degree in time and space. We use a software architecture that is built upon the C++ library \textit{deal.II} \cite{arndtDealIILibrary2020} combined with the linear algebra package \textit{Trilinos} \cite{thetrilinosprojectteamTrilinosProjectWebsite2020} under a Message Passing Interface (MPI) parallelization. All routines for the assembly of the Newton-linearized algebraic system and the linear algebraic solver are parallelized and can be run on standard multicore/-processor architectures. Here, the mesh is partitioned into the number of MPI processes and each of these processes owns the cells and corresponding matrix and vector entries of its partition (and some additional overhead). Basically, the simulation can be run with an arbitrary number of MPI processes. For the simulations of  this work we distributed exactly one MPI process to one physical CPU core.

Now we firstly address the integration over cut cells which is an important ingredient of the practical realization of the CutFEM approach. Then, key blocks of the code are discussed.  

\subsection{Integration over cut cells}
\label{Subsec:IterInt}

\begin{figure}[h!tb]
	\centering
    \subcaptionbox{\label{fig:kind_cuts_a}}
	[0.2\columnwidth]
    {\includegraphics[width=0.15\columnwidth,keepaspectratio]{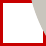}
    }
    \subcaptionbox{\label{fig:kind_cuts_b}}
	[0.2\columnwidth]
    {\includegraphics[width=0.15\columnwidth,keepaspectratio]{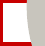}
    }
    \subcaptionbox{\label{fig:kind_cuts_c}}
	[0.2\columnwidth]
    {\includegraphics[width=0.15\columnwidth,keepaspectratio]{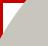}
    }
    \caption{Types of arising cut cells with gray shaded rigid body and white shaded fluid domain.}
    \label{fig:kind_cuts}
\end{figure}

In Problem~\ref{Prob_FullyDisc} or \ref{prob:Algeb_time_disc}, respectively, the stabilized discrete form $A_h^s$, that is defined in \eqref{Eq:AhS} along with \eqref{Def:Ah} and \eqref{Eq:Slf_a}, evokes the computation of integrals over the fluid domain. Due to the usage of unfitted meshes, integrals over the portions of the cut cells, that are filled with fluid, have thus to be evaluated. Fig.~\ref{fig:kind_cuts} illustrates schematically the types of cell intersections that can be induced by the motion of the rigid body. For quadrilateral finite elements, that we use in our implementation, the fluid cell portions can be pentagons, general quadrilaterals or triangles. We recall that we have assumed (cf.\ Assumption~\ref{assump:Triang}), that each face of the quadrilateral is cut at most once by the boundary of the rigid body (cf.\ Fig.~\ref{fig:kind_cuts}). The type of the cut cell is determined by evaluating the levelset function, for $t\in (0,T)$,
\begin{equation}
	\label{Eq:Levelset}
	\theta(\vec{x},t) 
	\begin{cases}
		< 0 & \forall \, \vec{x} \in \Omega_r^t \,, \\
		= 0 & \forall \, \vec{x} \in \partial \Omega_r^t \,,\\
		> 0 & \forall \, \vec{x} \in \Omega_f^t
	\end{cases}
\end{equation}
in the corners of each cell.

\begin{figure}[h!tb]
	\centering
    \subcaptionbox{\label{fig:cut_integration_1}}
	[0.2\columnwidth]
    {\includegraphics[width=0.15\columnwidth,keepaspectratio]{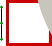}
    }
    \subcaptionbox{\label{fig:cut_integration_2}}
	[0.2\columnwidth]
    {\includegraphics[width=0.15\columnwidth,keepaspectratio]{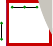}
    }
    %
    %
    \subcaptionbox{\label{fig:cut_integration_3}}
	[0.2\columnwidth]
    {\includegraphics[width=0.15\columnwidth,keepaspectratio]{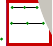}
    }
    \subcaptionbox{\label{fig:cut_integration_4}}
	[0.2\columnwidth]
    {\includegraphics[width=0.13\columnwidth,keepaspectratio]{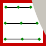}
    }
    \caption{Quadrature nodes of iterated numerical integration over cut cells by Gauss quadrature formulas.}%
    \label{fig:Cut_integration}
\end{figure}

In order to integrate over such cut cells we apply iterated one-dimensional numerical integration schemes. This is illustrated in Fig.~\ref{fig:Cut_integration} and reads as
\begin{equation}
\label{Def:ItQuad}
\iint_{K_{\text{fluid}}} f(x_1,x_2) \d (x_1,x_2) \approx \operatorname{meas}_2(K_{\text{fluid}}) \, \sum_{\mu_1=1}^{M_1} \sum_{\mu_2=1}^{M_2}  \omega_{\mu_1}  \omega_{\mu_2} f(x_{1,\mu_1},x_{2,\mu_2})\,,
\end{equation}
where $\omega_{\mu_1}$ and $x_{1,\mu_1}$, for $\mu_1=1,\ldots ,M_1$, as well as $\omega_{\mu_2}$ and $x_{2,\mu_2}$, for $\mu_2=1,\ldots ,M_2$, denote the weights and quadrature nodes of the quadrature formula for the respective coordinate direction and $K_{\text{fluid}}$ is the fluid portion of the cut cell $K$. Here, we use the Gauss quadrature formula. The number $M_1$ and $M_2$ of quadrature points in the coordinate directions, and thus the degree of exactness of the quadrature rules, can be chosen independently of each other, the number $M_2$ of quadrature nodes in $x_2$ direction can even depend on $\mu_1$, i.e., on the quadrature node $x_{1,\mu_1}$ in $x_1$ direction. We note that the integration by numerical quadrature is restricted to the fluid portions of the cut cells by adapting the interval lengths in either coordinate direction. The flexible choice of the degree of exactness of the iterated integration formulas then allows an accurate integration over cut cells. In the implementation, the direction $x_1$ and $x_2$ of the iterated integration is adapted to the respective element. Precisely, the direction of the largest element face (i.e., longest side in two space dimensions) bounding the fluid portion of the cut cell (cf.\ Fig.~\ref{fig:Cut_integration}), is chosen for the outer summation in Eq.~\eqref{Def:ItQuad}. For this we recall that in this work the assumption was made that the evolving domain and thus the position of the moving boundary $\Gamma_r^t$ (cf.\ Fig.~\ref{fig:problem_overview}) is prescribed explicitly and, therefore, is computable; cf.\ Assumption~\ref{assump:RB}. Further, a structured grid of elements aligned to the coordinate lines was supposed (cf.\ Subsec.~\ref{Subsec:SpaceDisc}). In our implementation, the quadrature points of each cut cell are computed in each time step of Problem~\ref{Prob_FullyDisc} before the assembly of the algebraic system is done. In the assembly routine, these points are then used to build a custom quadrature rule using \emph{deal.II} routines. 

Regarding the accuracy of the presented iterated integration we note the following. 
\begin{remark}
The numerical results presented in Sec.~\ref{sec:numerical_results} and \cite{anselmannNumericalConvergenceDiscrete2021} show that the space-time convergence behavior of the CutFEM approach is not deteriorated by the application of the iterated numerical integration on the cut cells. This also holds if irregular (tiny) cuts are present.  
\end{remark}

\begin{adjustbox}{valign=C,raise=\strutheight,minipage={1.0\linewidth}}
	\begin{wrapfigure}{r}{0.18\linewidth} 
		\centering    
		\includegraphics[width=1.0\linewidth]{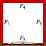}
		\caption{One face cut twice.}
		\label{fig:problem_curfed}
	\end{wrapfigure}%
	\strut{}
	\vspace*{-0.5cm} 
	\begin{remark}
		For rigid domains $\Omega_r^t$ with curved boundary, a cut scenario as sketched in Fig.~\ref{fig:problem_curfed} can occur. Here, the face $F_3$ of the cell is cut twice. If $\Omega_r$ is time-independent, such a scenario is avoided  by using the mesh sufficiently fine. Nevertheless, when the rigid domain $\Omega_r^t$ evolves in time, such cell intersections might arise. In our implementation, it is avoided by the following algorithm. We allow only the types of cell cuts that are illustrated in Fig.~\ref{fig:Cut_integration}.  To ensure the conformity of the arising cuts with these patterns, we evaluate the levelset function \eqref{Eq:Levelset} not only in the corners of each cell, but also in the $n$ Gauss-Lobatto points on each face of the corresponding cell. This is sketched in Fig.~\ref{fig:problem_curfed} for $n = 5$, which we use in all of the simulations presented in Sec.~\cref{sec:numerical_results}. If the levelset function coincides in both corners of the face, the values in the interior Gauss-Lobatto nodes of this face are checked to coincide with these value as well. If this condition is fulfilled, it is assumed that this face is not cut. Otherwise, the face is supposed to be cut more than once. For time-independent domains, either a mesh refinement or repositioning of the rigid domain $\Omega_r$ is performed. For evolving domains, a simple algorithm that slightly adjusts the time step size for this single time step in order to avoid such a scenario is applied.
		We explicitly note that cut scenarios as illustrated in Fig.~\ref{fig:problem_curfed} are very unlikely to occur, especially when the mesh size is chosen sufficiently fine.
	\end{remark}%
\end{adjustbox}
\begin{remark}
The integration over the boundary $\Gamma_r^t$ is performed using a standard parametrization of the circular rigid domain $\Omega_r^t$.
With this we transform the surface integral into a 1D integral that we numerically evaluate using a Gaussian quadrature rule.
\end{remark}

\subsection{Key code blocks}

\begin{figure}[h!tb] 
	\centering
	\includegraphics[width=1.0\textwidth,keepaspectratio]
	{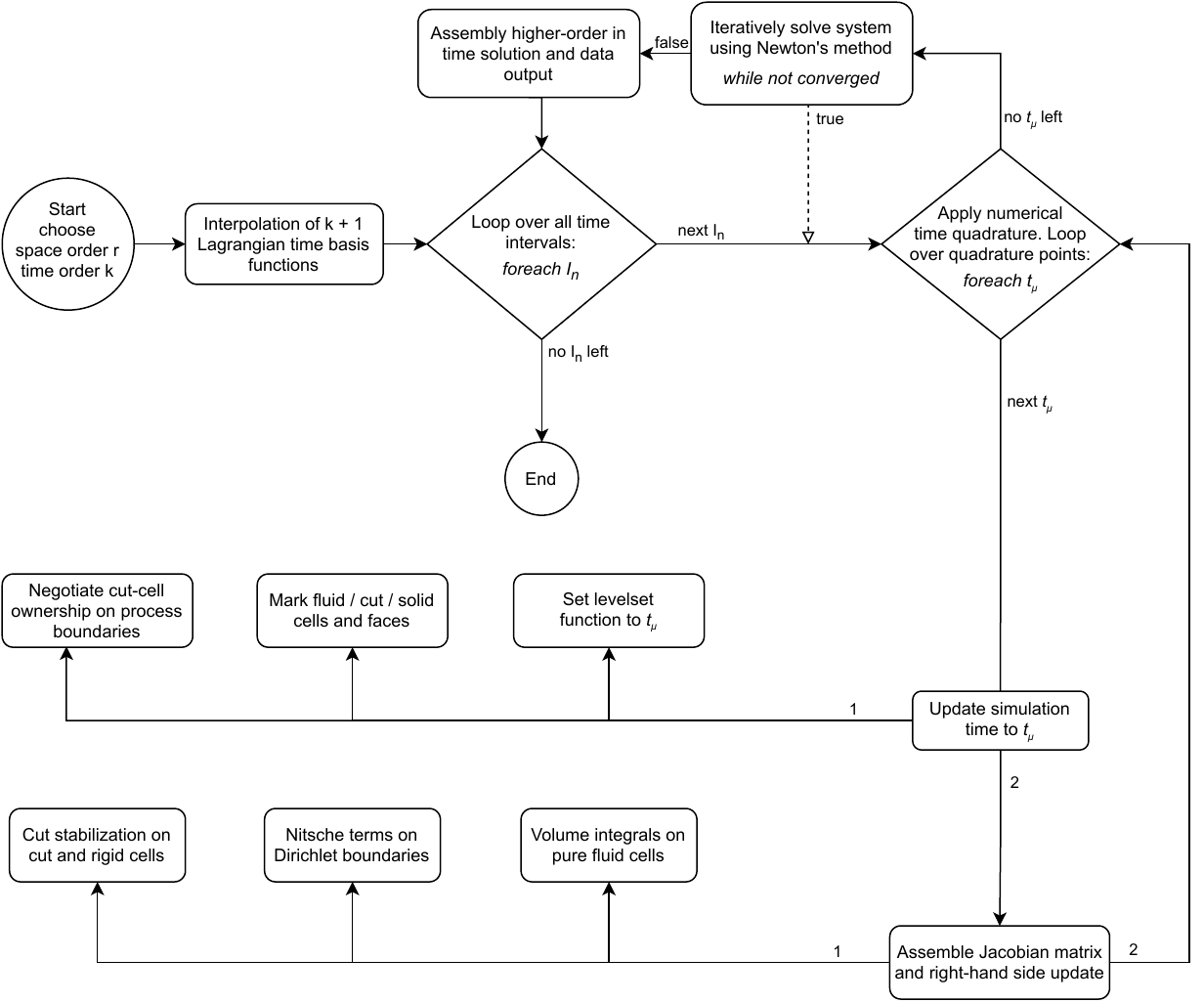}
	\caption{Key blocks of the CutFEM solver with arbitrary polynomial order in space and time.} 
	\label{fig:flow_chart}
\end{figure}

\begin{figure}[t]
	\centering
	\includegraphics[width=0.5\textwidth,keepaspectratio]
	{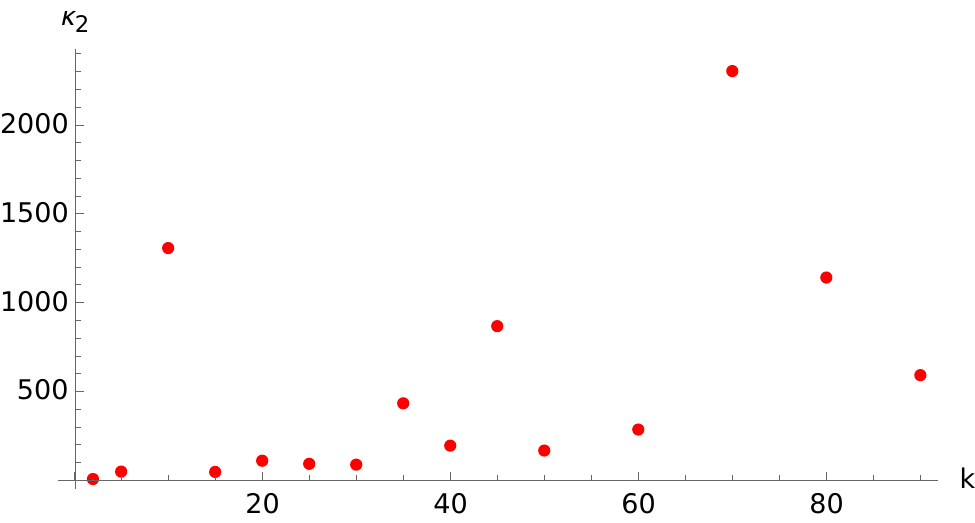}
	\caption{Spectral condition number $\kappa_2$ for the computation of the nodal $dG(k)$ time basis.} 
	\label{Fig:vandermonde_condition}
\end{figure}

Here, we describe the key blocks of the implementation of our CutFEM approach. A flowchart of the code is given in Fig.~\ref{fig:flow_chart}. After execution of the program and before the loop of the time marching process of Problem~\ref{Prob_FullyDisc} is started, a nodal Lagrangian time basis $\{\chi_l\}_{l=0}^k$ is computed with respect to the \mbox{$(k+1)$} Gauss--Radau quadrature nodes by solving the corresponding interpolation problem. Actually, this is done on the reference interval $\hat I = [0,1]$, such that $\{\hat \chi_l\}_{l=0}^k$, with $\hat \chi_l \in \mathbb P_k(\hat I;\R)$, satisfies
\[
\hat \chi_l(\hat t_{\mu}) = \delta_{l,\mu}\,, \quad l,\mu = 0\ldots, k\,,
\]
with the Kronecker symbol $\delta_{l,\mu}$ and the $(k+1)$ Gauss--Radau quadrature nodes $\{\hat t_{\mu}\}_{\mu=0}^l$ of $\hat I$. The resulting linear system is solved by using Linear Algebra PACKage (LAPACK) routines. Since the support points are non-equally distributed, the spectral condition number $\kappa_2$ of the resulting Vandermonde matrix is feasible even for high values of $k$. This is illustrated in Fig.~\ref{Fig:vandermonde_condition}.
In each subinterval of the time marching process of Problem~\ref{Prob_FullyDisc}, numerical quadrature is applied for the integration in the time and space domain. Firstly, the subinterval and its quadrature nodes are incremented to $I_n$ and $t_{n,\mu}$, for $\mu=0,\ldots,k$. For the representation of the rigid domain $\Omega_r^t$ we use the levelset function $\theta(\vec{x}, t)$, defined in \eqref{Eq:Levelset}. By evaluating this function we identify whether a finite element cell belongs at time $t_{n,\mu}$ to the set of fluid cells $\mathcal T_{h,f}^{t_{n,\mu}}$, rigid cells $\mathcal T_{h,r}^{t_{n,\mu}}$ or cut cells $\mathcal T_{h,c}^{t_{n,\mu}}$. In subroutines, we mark each cell of the computational background grid by evaluating the levelset function in the edges of each cell. For cut cells, we pre-compute the quadrature points and weights for the (fluid) volume and surface integrals (line integrals in two space dimensions) and store them into lookup tables.
In the MPI based implementation these tasks scale perfectly with the number of available processes.
\begin{remark}
	We note that the levelset function \eqref{Eq:Levelset} is not discretized or computed by the approximation of a suitable transport equation. The direct evaluation of \eqref{Eq:Levelset} becomes feasible due to the Assumption~\ref{assump:RB} that the motion of the rigid body is precribed. Thereby, further discretization errors are avoided. 
\end{remark}

\begin{figure}[h!tb]
	\centering
	\includegraphics[width=0.5\textwidth,keepaspectratio]
	{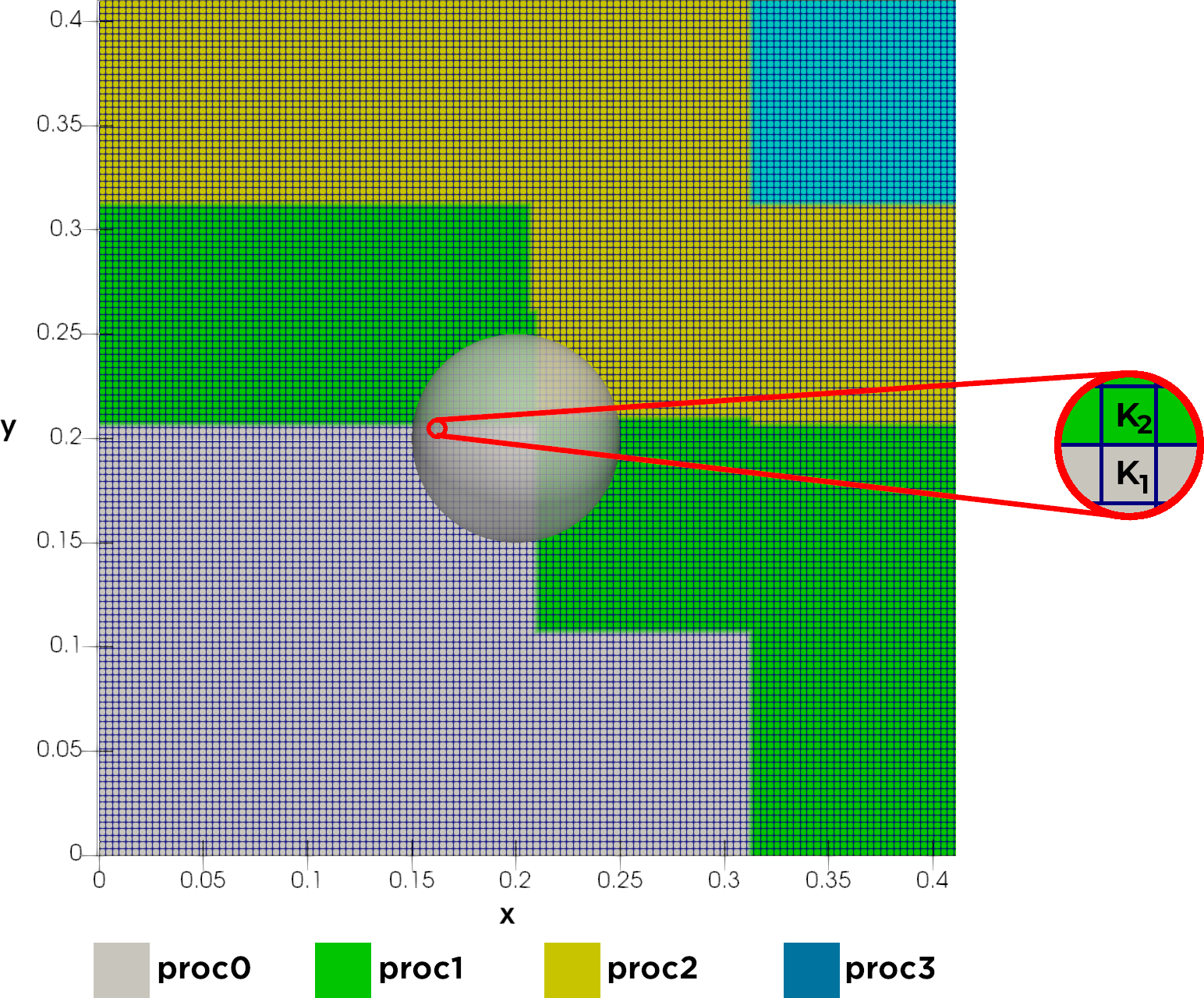}
    \caption{Assembly of the ghost penalty stabilization at partition boundaries with computational mesh $\mathcal{T}_h$ and colored processor ownerships of mesh cells.} 
	\label{Fig:process_marking}
\end{figure}

However, the parallelized code needs some mechanism to prevent that the ghost penalty operator defined in Eq.~\eqref{eq:extension_operator} is applied twice in the interface zone of MPI process boundaries of the mesh partition. For this we use a simple master-slave approach and let the process with the higher process number assemble the ghost penalty operator over patches $\omega_F$ at process partition boundaries. This problem is sketched in Fig.~\ref{Fig:process_marking} for the patch $\omega_F$, built from the cells $K_1$ and $K_2$. Cell $K_1$ of the patch belongs to process 0 and cell 2 to process 1. Without a control of the parallel assembly of the ghost penalty stabilization, both processes would assemble the contributions of Eq.~\eqref{eq:extension_operator}, such that the stabilization would be applied twice. To avoid this, the process with the higher process number (process 1 in this case) assembles the ghost penalty stabilization.

In the next block the linear system of the Newton iteration (cf.\ Remark~\ref{Rem:AlgSolv}) is assembled. For the integration over fluid cut cells iterated integration, as sketched in Subsec.~\ref{Subsec:IterInt}, is applied. Then, the resulting linear system is solved by the parallel, sparse direct solver {SuperLU\_DIST} \cite{liSuperLUDISTScalable2003} and, unless the stopping criteria is satisfied, the Newton iteration is continued by reassembling the linear system. After convergence of the Newton iteration, a parallel assembly and data output of the computed higher order solution for the subinterval $I_n$ terminates the time step.

\section{Numerical results}
\label{sec:numerical_results}

In this section we present the results of four numerical experiments to illustrate the performance and accuracy properties of the proposed CutFEM approach. We start with a space-time convergence test with a moving $\Omega_r^t$, to validate all aspects of the implementation.
In the second experiment, the popular DFG benchmark of flow around a cylinder (cf.~\cite{schaferBenchmarkComputationsLaminar1996}) is considered.
A careful comparative study of the CutFEM approach on background meshes with the results obtained for body fitted meshes is done.
In the last two test cases we investigate the quality of the ghost penalty stabilization and extension of \eqref{eq:extension_operator} on an evolving domain and demonstrate it's application in a realistic example. Throughout all simulations we set the viscosity to $\nu = 0.001$.

In all the examples, the Newton solver uses an absolute value of the residual smaller than $\num{1e-10}$ as the stopping criteria. In \cref{sec:flow_re_100}, it is checked additionally if the initial residual of a time step is decreased by a factor of 1000. For all our numerical experiments we set the Nitsche penalty parameters of \Cref{Def:B_GamD} to $\gamma_1 = \gamma_2 = 35$, which is motivated by the numerical experiments in \cite{schottNewFaceorientedStabilized2014,winterNitscheCutFinite2018}.
The ghost penalty parameters of \Cref{eq:extension_operator} are set to $\gamma_v = \gamma_p = 10^{-2}$, except in \cref{sec:dynamic_flow_around_moving_cylinder}, where we set $\gamma_v = \gamma_p = 1$.

\subsection{Experimental order of convergence}

\begin{figure}[h!t]
	\centering
	\subcaptionbox{Initial configuration of the space-time convergence test. \label{fig:experiment_setup}}
	[0.49\columnwidth]
	{\includegraphics[width=0.3\textwidth,keepaspectratio]{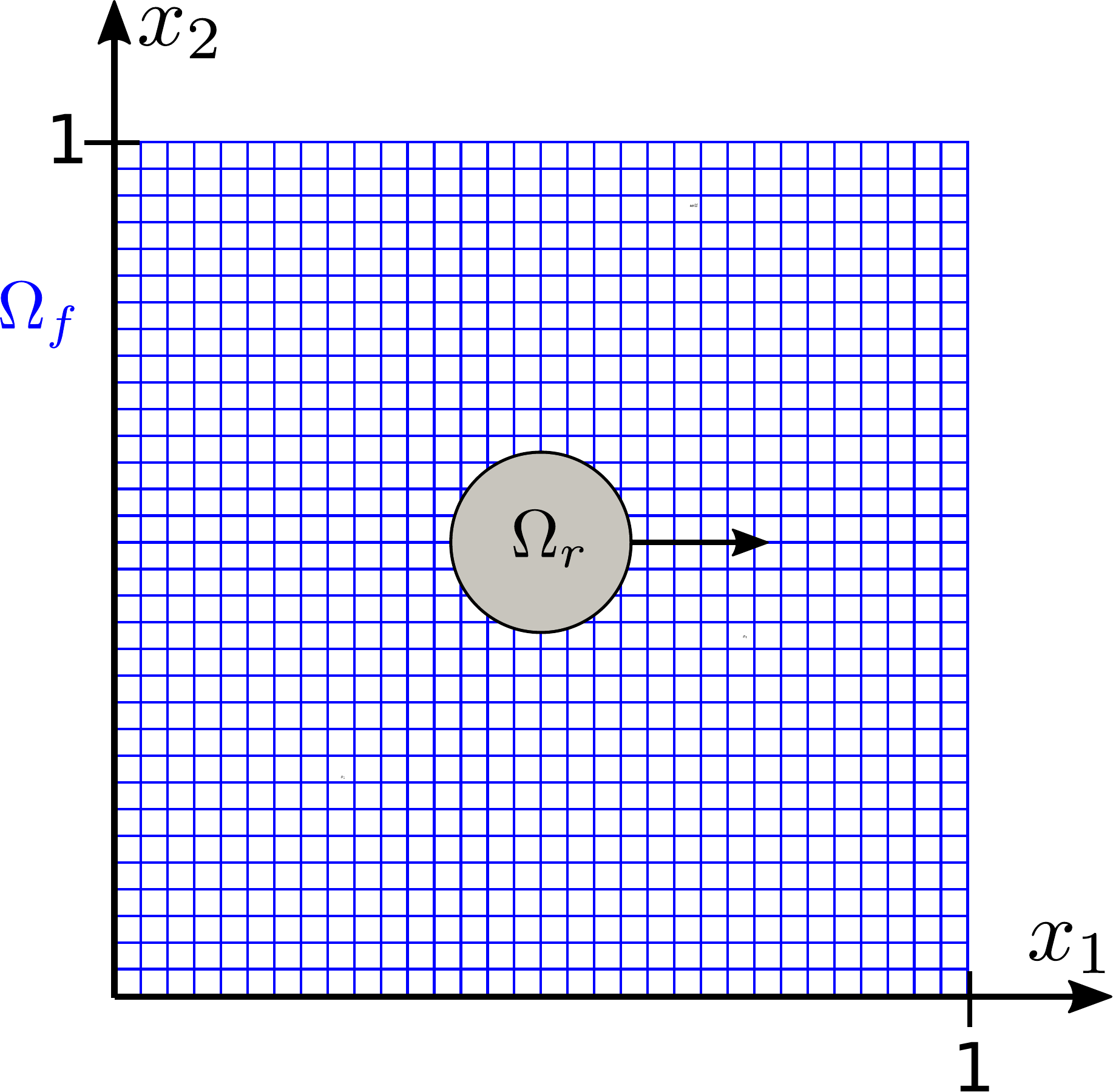}
	}
	\subcaptionbox{Solution of the pressure field at $t=T$ with $h = \frac{h_0}{2^4}$. \label{fig:experiment_setup_final}}
	[0.49\columnwidth]
	{\includegraphics[width=0.3\textwidth,keepaspectratio]{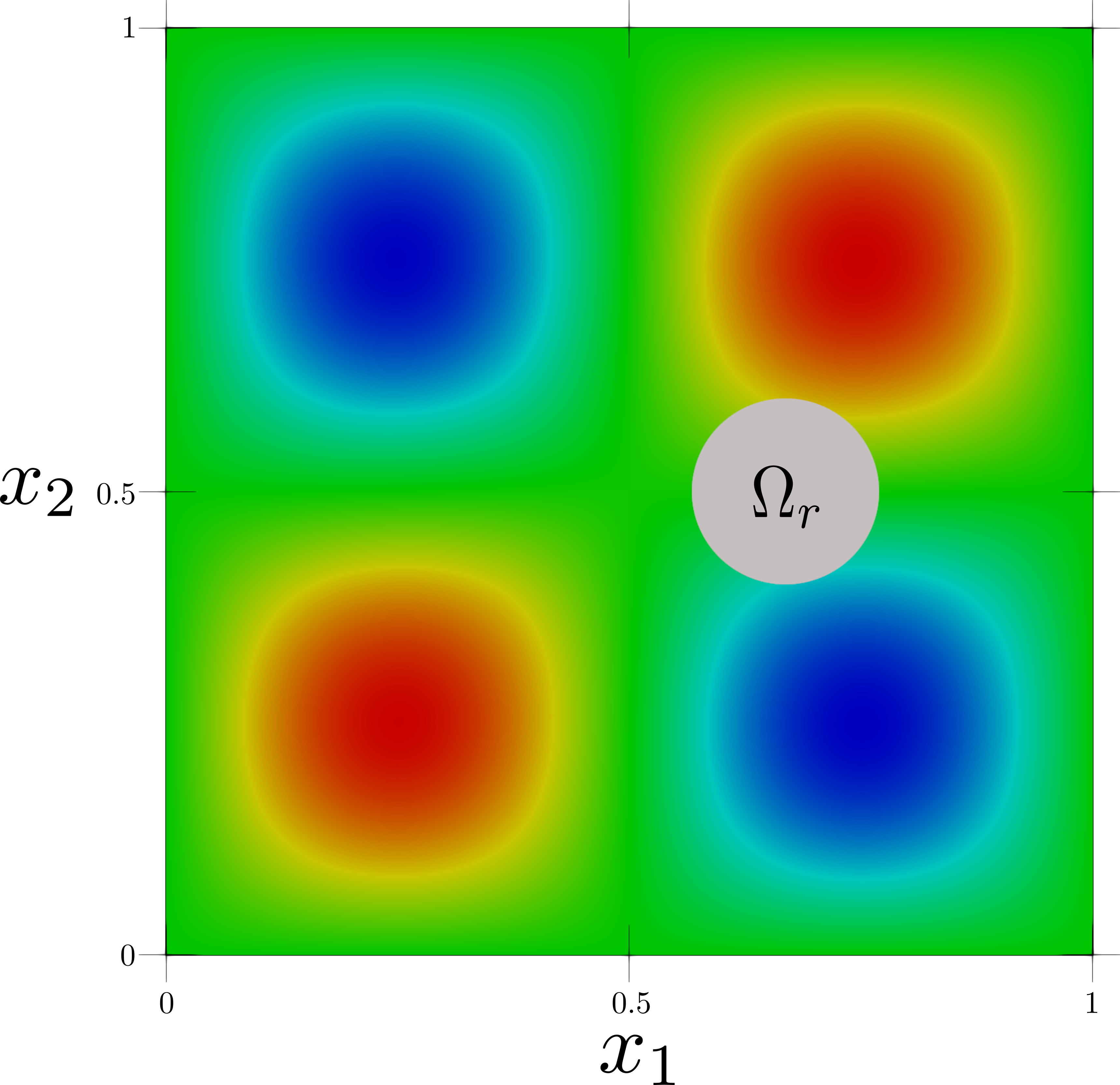}
	}
	\caption{Problem setup and pressure field at the end of the simulation.}
	\label{fig:moving_convergence_test}
\end{figure}

The aim of this first numerical experiment is to verify all components of our approach, in particular the integration over cut cells as well as the ghost penalty stabilization, by a numerical space-time convergence study. We consider the problem setting sketched in Fig.~\ref{fig:experiment_setup} with a time-dependent domain $\Omega_r^t$. We put $\Omega \times I = (0,1)^2 \times (0,1]$.
The rigid body is modeled by a circular domain of radius of $r = 0.1$ and an initial position of its center at $\vec{x}_r(0) = (0.5, 0.5)^\top$.
The motion of the center is prescribed by 
\begin{equation}
\label{eq:rigid_motion}
	\vec{x}_r(t) = \vec{x}_r(0) + \left( A \cdot \sin(\omega \cdot t), 0 \right)^\top\,,
\end{equation}
with $A = 0.2$ and $\omega = 1$. We prescribe the initial value $\vec{v}_0$ and right-hand side function $\vec{f}$ on $\Omega_f \times I$ in such a way, that the solution of the Navier--Stokes system on $\Omega_{f}$ is given by
\begin{align*}
	\vec{v}_e(\vec{x},t)
	&\coloneqq
	\begin{pmatrix}
		\cos(x_2 \pi) \cdot \sin(t) \cdot \sin(x_1 \pi)^2  \cdot \sin(x_2 \pi) \\
		- \cos(x_1 \pi) \cdot \sin(t) \cdot \sin(x_2 \pi)^2  \cdot \sin(x_1 \pi)
	\end{pmatrix} \,,
	\\
	p_e(\vec{x},t)
	&\coloneqq
	\cos(x_2 \pi) \cdot \sin(t) \cdot \sin(x_1 \pi) \cdot \cos(x_1 \pi) \cdot \sin(x_2 \pi) \,.
\end{align*}
On the inner fluid boundary we prescribe the condition $\vec{g}_r = \vec{v}_e$. On the outer boundary we use a homogeneous Dirichlet condition.
Table~\ref{tab:numerical_results_1} shows the computed errors and experimental orders of convergence (EOC) for two different combinations of space-time finite elements, based on the Taylor--Hood family in space. In the first experiment, the polynomial order in time is adapted to the spatial apporximation of the prerssure variable. It is chosen to $k=1$ such that optimal second order of convergence in time and space is obtained for the pressure variable with discrete values in $H_h^1$. In the second experiment, the polynomial order in time is put to $k=2$ such that optimal third order of convergence in time and space is obtained for the pressure variable with discrete values in $H_h^2\times H_h^2$. For this we recall that a non-equal order, inf-sup stable discretization in space by the Taylor--Hood family of elements is used here. For the definition of the discrete function spaces we refer to \eqref{Eq:DefHh} and \eqref{Def:PkIB}, respectively. In Table~\ref{fig:moving_convergence_test}, the expected convergence of optimal order is documented for both experiments. 

\sisetup{scientific-notation = true,
	round-mode=places,
	round-precision=3,
	output-exponent-marker=\ensuremath{\mathrm{e}},
	table-figures-integer=1, 
	table-figures-decimal=3, 
	table-figures-exponent=1, 
	table-sign-mantissa = false, 
	table-sign-exponent = true, 
	table-number-alignment=center} 
\begin{table}[!h]
	\caption{Errors and experimental order of convergence for $\tau_0 = 1.0$ and $h_0 = {1}/({2 \sqrt{2}})$.}
	\centering
	\begin{tabular}{c@{\,\,\,\,}c  S@{\,}c  S@{\,}c  S@{\,}c  S@{\,}c}
		\toprule
		{$\tau$} & {$h$} &
		{ $\| e^{\;\vec{v}}  \|_{L^2(L^2)} $ } & {EOC} &
		{ $\| e^{\;p}  \|_{L^2(L^2)} $ } & {EOC} &  
		{ $\| e^{\;\vec{v}}  \|_{L^2(L^2)} $ } & {EOC} &
		{ $\| e^{\;p}  \|_{L^2(L^2)} $ } & {EOC} \\
		\cmidrule(r){1-2}
		\cmidrule(lr){3-4}
		\cmidrule(lr){5-6}
		\cmidrule(l){7-8}
		\cmidrule(l){9-10}
		$\tau_0/2^0$ & $h_0/2^0$ & 4.3136e-02 & {--} & 1.4756e-02 & {--} & 3.4198e-02 & {--} & 1.3397e-02 & {--}  \\ 
		$\tau_0/2^1$ & $h_0/2^1$ & 8.6019e-03 & 2.33 & 3.0469e-03 & 2.28 & 7.8334e-03 & 2.13 & 2.5895e-03 & 2.37  \\
		$\tau_0/2^2$ & $h_0/2^2$ & 2.0628e-03 & 2.06 & 6.9638e-04 & 2.13 & 9.4588e-04 & 3.05 & 5.6048e-04 & 2.21  \\
		$\tau_0/2^3$ & $h_0/2^3$ & 4.8824e-03 & 2.08 & 1.8527e-04 & 1.91 & 1.1506e-04 & 3.04 & 1.3538e-04 & 2.05  \\
		$\tau_0/2^4$ & $h_0/2^4$ & 1.2206e-04 & 2.00 & 4.5999e-05 & 2.01 & 1.4283e-05 & 3.01 & 3.3846e-05 & 2.00  \\
		\cmidrule(r){1-2}
		\cmidrule(lr){3-6}
		\cmidrule(lr){7-10}
		\multicolumn{2}{c}{$\text{elements}_{|I_n}$} &
		\multicolumn{4}{c}{$(\mathbb P_1(I_n;H_h^2))^2 \times \mathbb P_1(I_n;H_h^1)$} &
		\multicolumn{4}{c}{$(\mathbb P_2(I_n;H_h^2))^2 \times \mathbb P_2(I_n;H_h^1)$} \\
		\bottomrule
	\end{tabular}
	\label{tab:numerical_results_1}
\end{table}
\begin{remark}
	The setting of this example was chosen in such a way, that the propagation of the boundary of $\Gamma_r^t$ is less than one cell in each time step, in order to preserve a CFL like condition. We didn't observe any decrease in the convergence rate by violating this condition. For instance, choosing $h_0 = 1/(4\sqrt{2})$, leads to comparable EOCs, but violates the assumption that the boundary of $\Gamma_r^t$ propagtes over a disctance less than one cell within a time step.
\end{remark}

\subsection{Time periodic flow around a cylinder}
\label{sec:flow_re_100}
The aim of the second numerical example is to analyze the effects of usage of cut cells for the space discretization on background meshes along with the extension to a rigid domain $\Omega_r$ and the application of the stabilization introduced in \eqref{eq:extension_operator}.
For this, we use the well-known DFG benchmark setting of flow around a cylinder, defined in \cite{schaferBenchmarkComputationsLaminar1996}.
Although the domain is non-evolving, evaluating the performance properties of the approach for this flow benchmark is of high interest. We compare the numerical results with the ones obtained for simulations on a body-fitted meshes that is highly pre-adapted to the cylinder. Quantities of interest and comparison in the simulations are the drag and lift coefficient of the flow on the circular cross section (cf.~\cite{schaferBenchmarkComputationsLaminar1996}). With the drag and lift forces $F_D$ and $F_L$ on the rigid circle $S$ given by 
\begin{align}
	F_D &= \int_S
	\left(
	\nu \frac{\partial v_t}{\partial \vec{n}} n_y - P n_x
	\right) \d S \,,
	&
	F_L &= - \int_S
	\left(
	\nu \frac{\partial v_t}{\partial \vec{n}} n_x - P n_y
	\right) \d S \,,
\end{align}
where $\vec n$ is the normal vector on $S$ and $v_t$ is the tangential velocity $\vec{t} = (n_y, -n_x)^\top$, the drag and lift coefficient $c_D, c_L$ are defined by means of
\begin{align}
	c_D &= \frac{2}{\bar{U}^2 L} F_D \,,
	&
	c_L &= \frac{2}{\bar{U}^2 L} F_L \,.
\end{align}

According to \cite{schaferBenchmarkComputationsLaminar1996}, we set the boundary condition on the inflow boundary $\Gamma_i$ as $\vec{g}_i(x, y, t) = (\frac{4 \cdot 1.5 \cdot y (0.41 - y)}{0.41^2}, 0)^\top$. This leads to a Reynold's number of $Re = 100$ and a time-periodic flow behavior. The space--time discretization is done in the discrete spaces $\mathbb P_1(I_n;H_h^2))^2 \times \mathbb P_1(I_n;H_h^1)$ on each time interval $I_n$. For the CutFEM approach we study a sequence of successive mesh refinements in space. For the developed flow and the computation of the drag and lift coefficient as quantities of physical interest we use the time step size $\tau = \num{0.005}$.
\begin{figure}[h!t]
	\centering
    \subcaptionbox{Cut$_7$ configuration with computational mesh. \label{fig:cut_grid}}
	[0.49\columnwidth]
	{\includegraphics[width=0.4\textwidth,keepaspectratio]{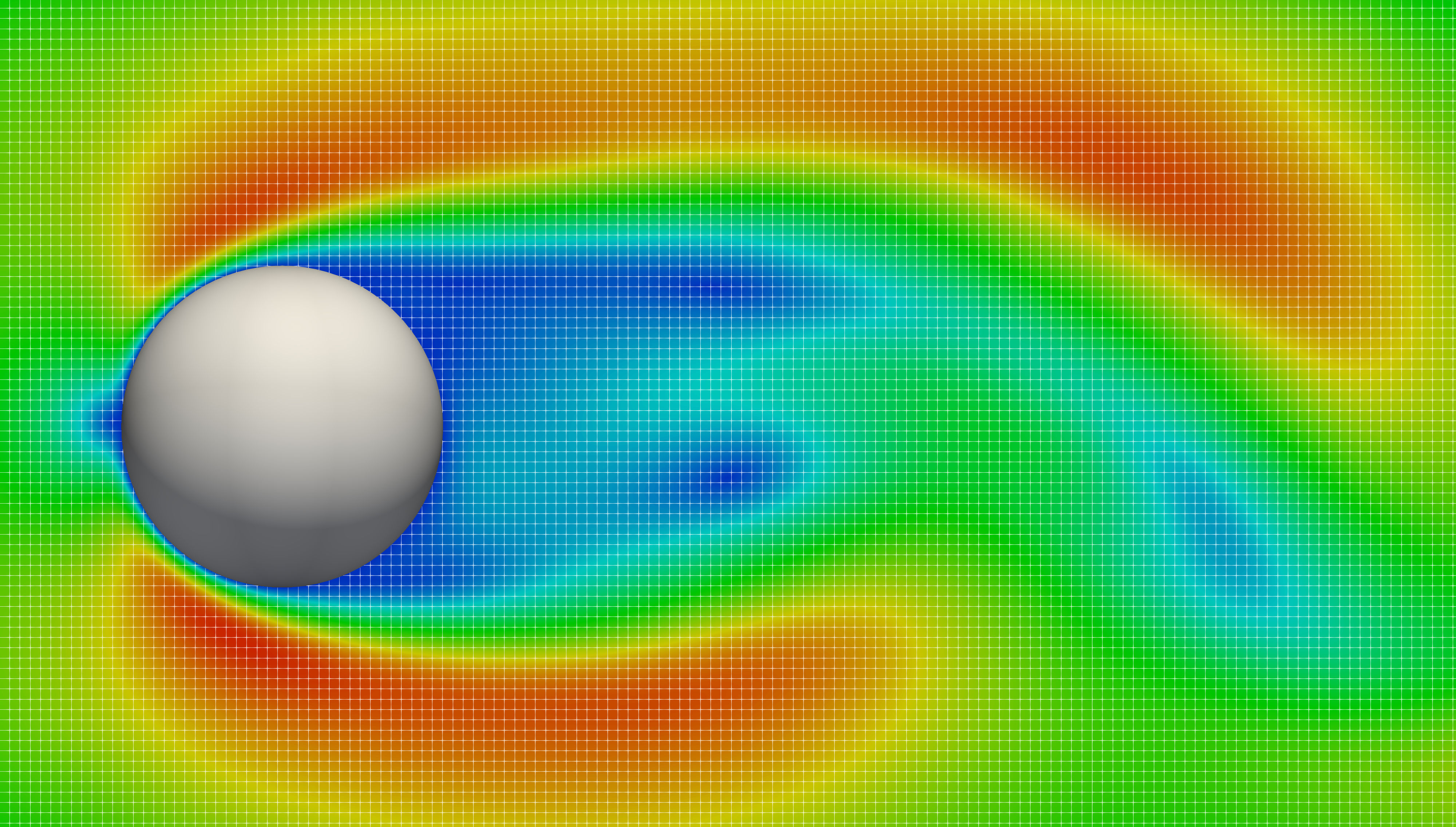}
	}
    \subcaptionbox{Fitted mesh, pre-adapted to the rigid body.\label{fig:fitted_grid}}
	[0.49\columnwidth]
    {\includegraphics[width=0.4\textwidth,keepaspectratio]{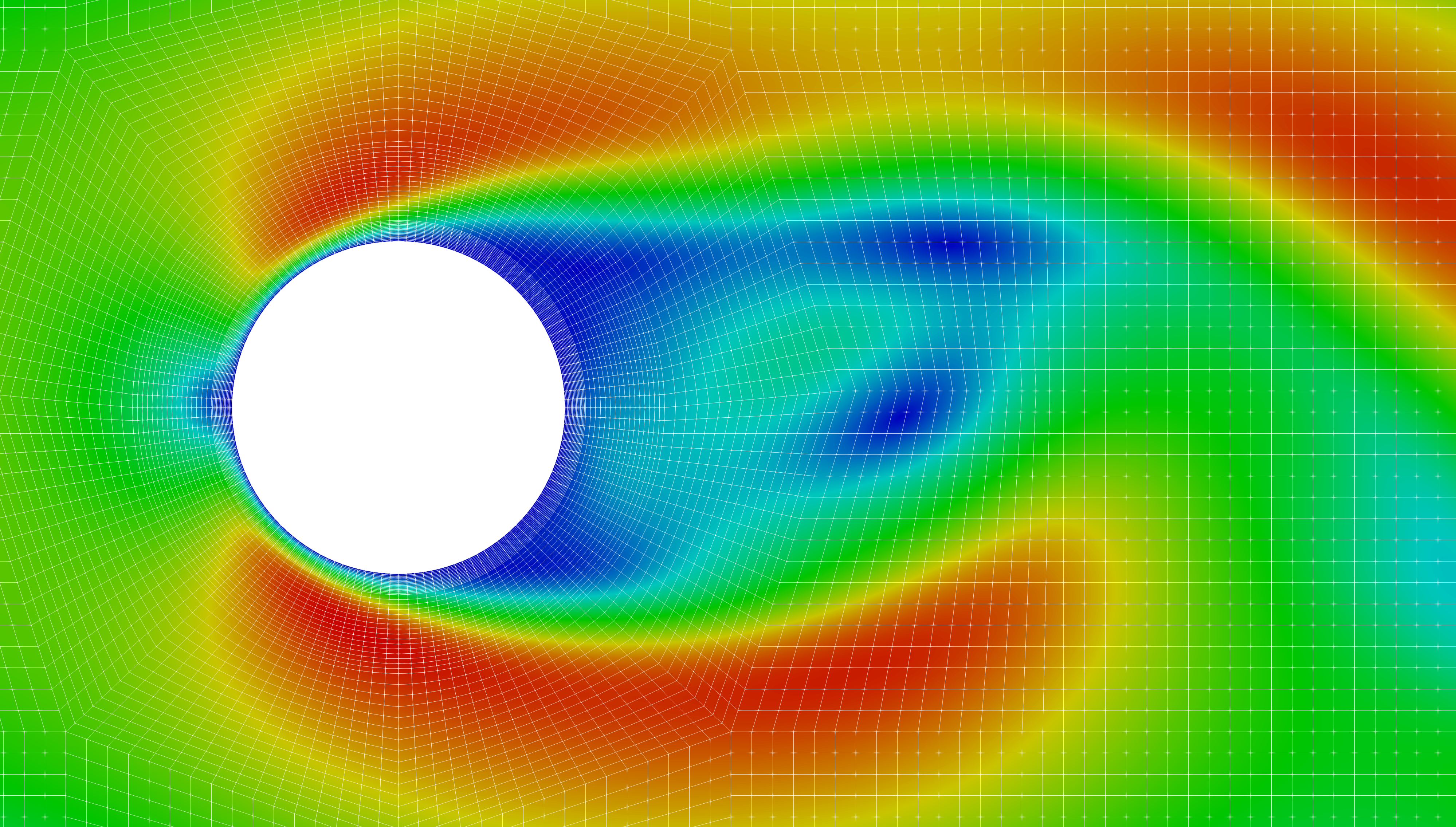}
	}
    \caption{Portion of problem setting and computed velocity profile for CutFEM approach (\Cref{fig:cut_grid}) and fitted mesh approximation (\Cref{fig:fitted_grid}).}
	\label{fig:dfg_benchmark_grids}
\end{figure}

\begin{figure}[h!t]
	\centering
	\subcaptionbox{Drag coefficients $c_D$ \label{fig:drag_coefficients}}
	[1.\columnwidth]
	{\includegraphics[width=0.7\textwidth,keepaspectratio]{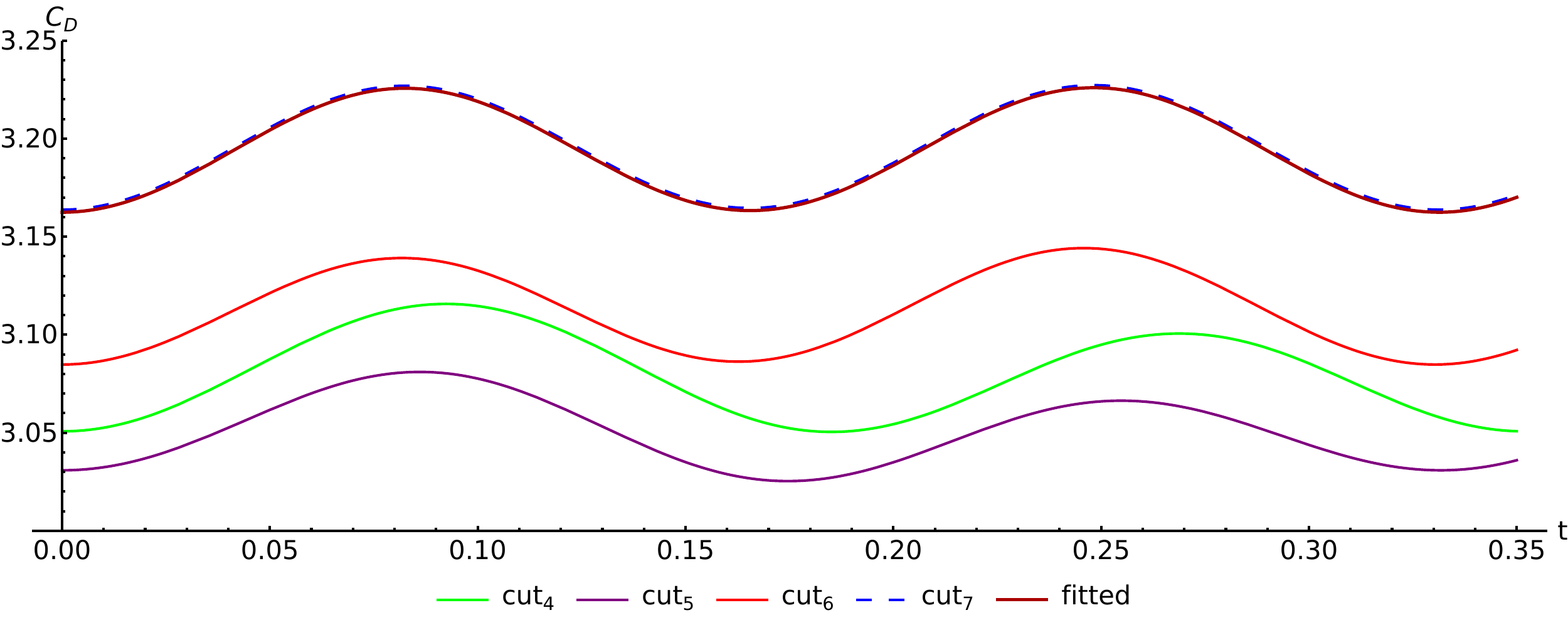}
	}\\[5ex]
	\subcaptionbox{Lift coefficients $c_L$ \label{fig:lift_coefficients}}
	[1.\columnwidth]
	{\includegraphics[width=0.7\textwidth,keepaspectratio]{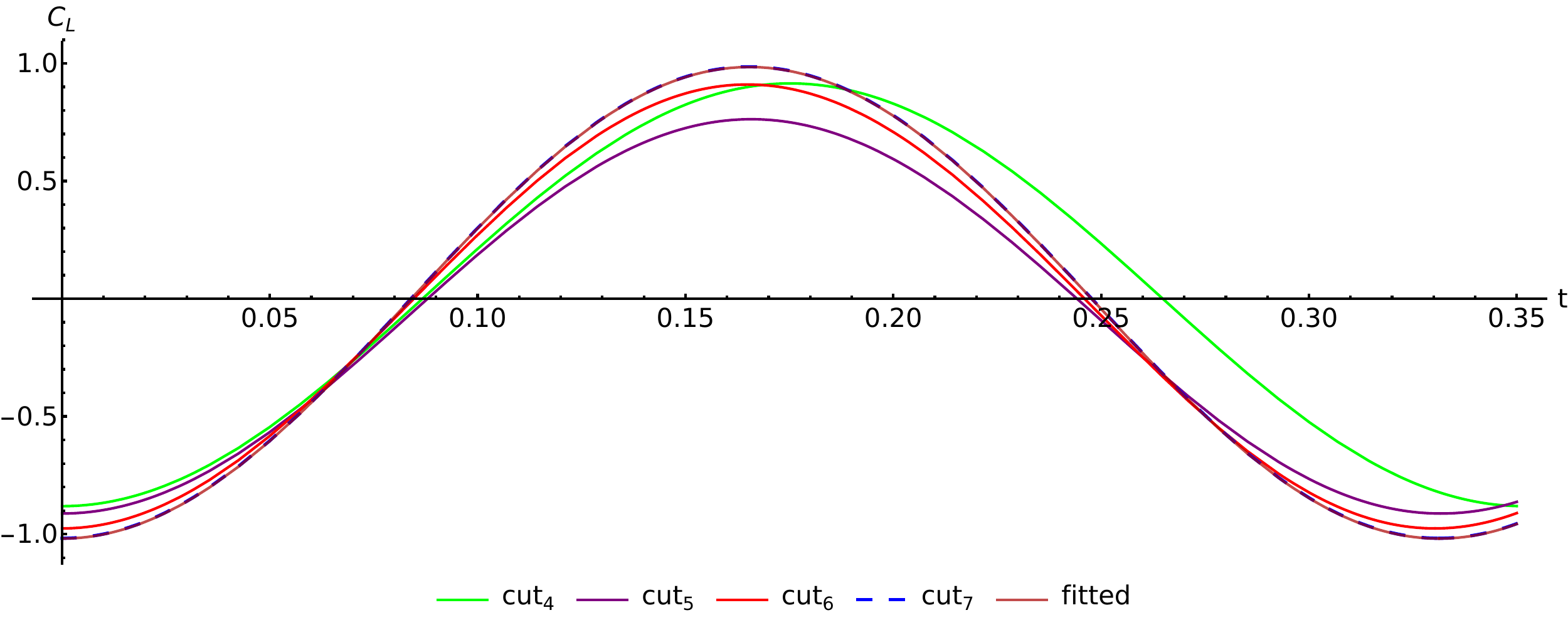}
	}
    \caption{Computed drag and lift coefficients for the example of \Cref{sec:flow_re_100} and different time discretization schemes with basic time step size $\tau = \num{0.005}$.}
	\label{fig:drag_lift}
\end{figure}

Fig.~\ref{fig:dfg_benchmark_grids} shows the computed flow profile and spatial mesh of a CutFEM simulation and a computation done on a highly pre-adapted, body-fitted mesh. In the latter case, high accuracy for the drag and lift coefficient is obtained (cf.~\cite{schaferBenchmarkComputationsLaminar1996}). The CutFEM simulation is done for a sequence of successive refinement steps of the entire background mesh (cf.\ Table~\ref{tab:drag_lift_coefficients}). Table~\ref{tab:drag_lift_coefficients} shows the computed drag and lift coefficient as well as their frequency. Precisely, the latter denotes the  frequency of the oscillation of the lift coefficient $c_L$, which is computed at a time instant $t_0$ where the lift coefficient $c_L$ is smallest and ends at a time instant $t_1 = t_0 + \tfrac{1}{f}$ when $c_L$ is smallest again. Exemplarily, Fig.~\ref{fig:drag_lift} illustrates such a cycle for the drag and the lift coefficient. To simplify the interpretation of the graphs, the starting point of the monitoring cycle is set to $t = 0$.
{
\sisetup{scientific-notation = false,
	     round-mode=places,
         round-precision=4,
		 table-figures-integer=1, 
		 table-figures-decimal=4, 
	 	 table-number-alignment=center-decimal-marker} 
\begin{table}[h!t]
    \small
    \caption{Computed minimum and maximum drag--lift coefficients and frequency for a sequence of successively refined CutFEM background meshes and an adapted, body-fitted mesh.}
	\centering
	\begin{tabular}{ccSSSSS}
		\toprule
        Method & $\text{DoFs}_{|I_n}$ & {min($c_L$)} & {max($c_L$)} & {min($c_D$)} & {max($c_D$)} & {$f_L$} \\
        \cmidrule(r){1-2}
        \cmidrule(lr){3-4}
        \cmidrule(lr){5-6}
        \cmidrule(l){7-7}
        cut$_4$ & \num{24000} & -0.881383 & 0.914682 & 3.05072 & 3.11567 & 2.84333 \\
        cut$_5$ & \num{94086} & -0.912680 & 0.762557 & 3.03084 & 3.08099 & 3.01459 \\
        cut$_6$ & \num{372486} & -0.975867 & 0.910129 & 3.08479 & 3.13906 & 3.02737 \\
        cut$_7$ & \num{1482246} & -1.01765 & 0.985596 & 3.16375 & 3.22723 & 3.01825 \\[2ex]
        \begin{tabular}{@{}l@{}} Adapted,\\ body-fitted \end{tabular} & \num{502464} & -1.01899 & 0.984292 & 3.16245 & 3.22593 & 3.01844 \\
		\bottomrule
	\end{tabular}
\label{tab:drag_lift_coefficients}
\end{table}
}

Comparing the results of the either approaches, CutFEM on background meshes versus adapted body-fitted mesh, we observe that the goal quantities of the CutFEM simulations nicely converge to the results computed on a manually pre-adapted and body-fitted mesh. Whereas the CutFEM background mesh is refined uniformly, the body-fitted mesh is highly adapted to the problem setting and flow profile. In particular, the region around the cylinder is strongly resolved by the spatial mesh whereas the CutFEM mesh remains much coarser in this region, even on the finest refinement level. 
In Fig.~\ref{fig:dfg_benchmark_grids} both meshes are shown. The diameter of a quadratic cut cell of the $\text{Cut}_7$ configuration is \num{0.0045299}. The diameter of a rectangular cell of the body-fitted triangulation in the neighborhood of $\Omega_r$ is \num{0.00249503}.
Therefore, the superiority of the body-fitted approach for this test case with non-evolving domain is obvious. Nevertheless, Table~\ref{tab:drag_lift_coefficients} and Fig.~\ref{fig:drag_lift} nicely show the numerical convergence of the CutFEM simulations for a successive refinement of the background mesh to the expected solution in terms of the confirmed drag and lift coefficient and their frequency (cf.~\cite{schaferBenchmarkComputationsLaminar1996}). This confirms the accuracy of the proposed CutFEM approach in the case of non-evolving domains.

\subsection{Dynamic Poiseuille flow}

In this numerical experiment the accuracy of the presented CutFEM approach is analyzed for an evolving domain.
Precisely, we study the perturbation of a Poiseuille flow profile by the motion of an enclosed moving rigid body (cf.\ Fig.~\ref{fig:problem_overview}) and whether application of the ghost penalty operator leads to a disturbance of the solution.
The test configuration is chosen in such a way, that the Poiseuille profile is preserved in the pipe, even though a moving rigid body is present. We than analyze if the Poiseuille flow is recovered by the CutFEM approach. 

The computational domain is $\Omega = (0, l_x) \times (-l_y, l_y)$, with $l_x = 3, l_y = 0.5$, we let $I = (0, 40]$. A parabolic Poiseuille flow profile is prescribed at the left inflow boundary, that is given by  
\begin{equation}
	\label{eq:PoisFlow}
	\vec{v}_p(x,y) = U_{\text{in}}(l_y^2 - y^2, 0)^\top
\end{equation}
with $U_{\text{in}} = 1$. The initial value is $\vec v_0 = \vec 0$. The extension of \eqref{eq:PoisFlow} to the entire pipe is the Poiseuille flow and satisfies the Navier--Stokes equations along with the pressure $p_p = - 2 \nu U_{\text{in}} (x - l_x) $. The rigid body is modeled by a circular domain of radius of $r = 0.2$ and an initial position of its center at $\vec{x}_r(0) = (1.545, 0)^\top$.
The motion of the ball's center is prescribed by \cref{eq:rigid_motion}, with $A = 0.8$ and $\omega = 0.2$. On the boundary of the rigid body the Dirichlet condition   
\begin{equation}
	\label{eq:BC_RB}
	\vec v = \vec{g}_r := \vec{v}_p
\end{equation}
is prescribed. After a transition from the zero initial state, the Poiseuille flow profile is expected to develop, and it is not perturbed by the rigid body motion due to the choice of the boundary condition \eqref{eq:BC_RB}. Clearly, even though the rigid body is moving inside the fluid by means of \eqref{eq:rigid_motion}, its impact on the fluid flow is hidden by means of the boundary condition \eqref{eq:BC_RB}. 

\begin{figure}[tb]
    \centering
    \includegraphics[width=0.6\linewidth]{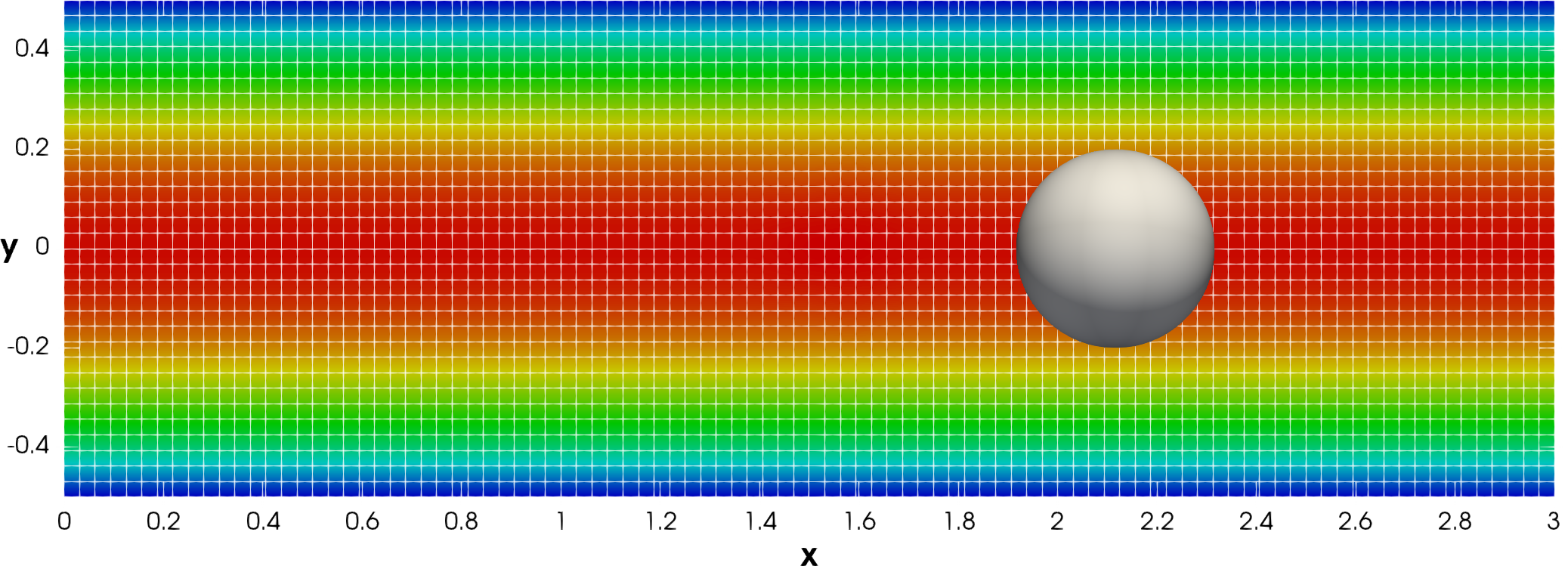}
    \caption{Problem setting and flow profile at $T = 40$}%
    \label{fig:poiseuille_setting}
\end{figure}
\begin{figure}[!htb]
	\centering
    \subcaptionbox{Velocity along $y$-direction at $x = 2.34$ \label{fig:poiseuille_velocity}}
	[1.\columnwidth]
    {\includegraphics[width=0.6\textwidth,keepaspectratio]
    {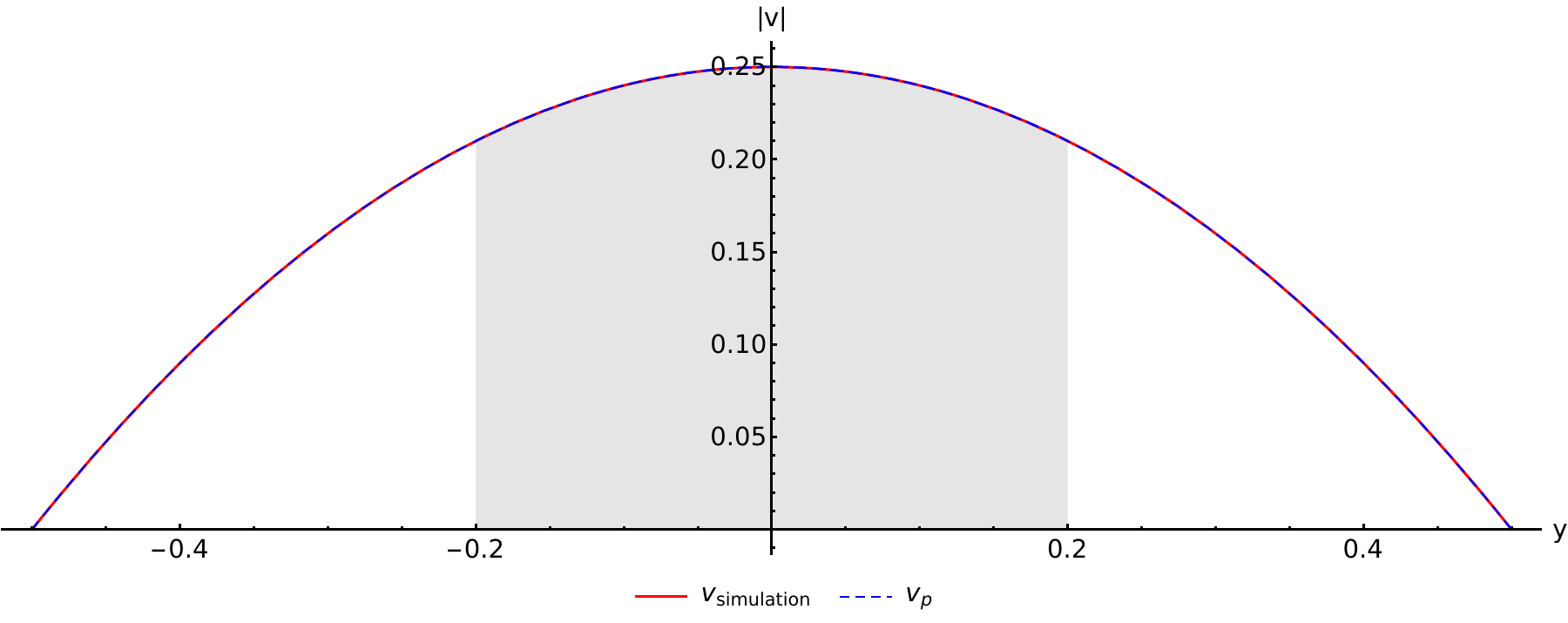}}
    %
    %
    \subcaptionbox{Pressure along $x$-direction at $y = 0$  \label{fig:poiseuille_pressure}}
	[1.\columnwidth]
    {\includegraphics[width=0.6\textwidth,keepaspectratio]
    {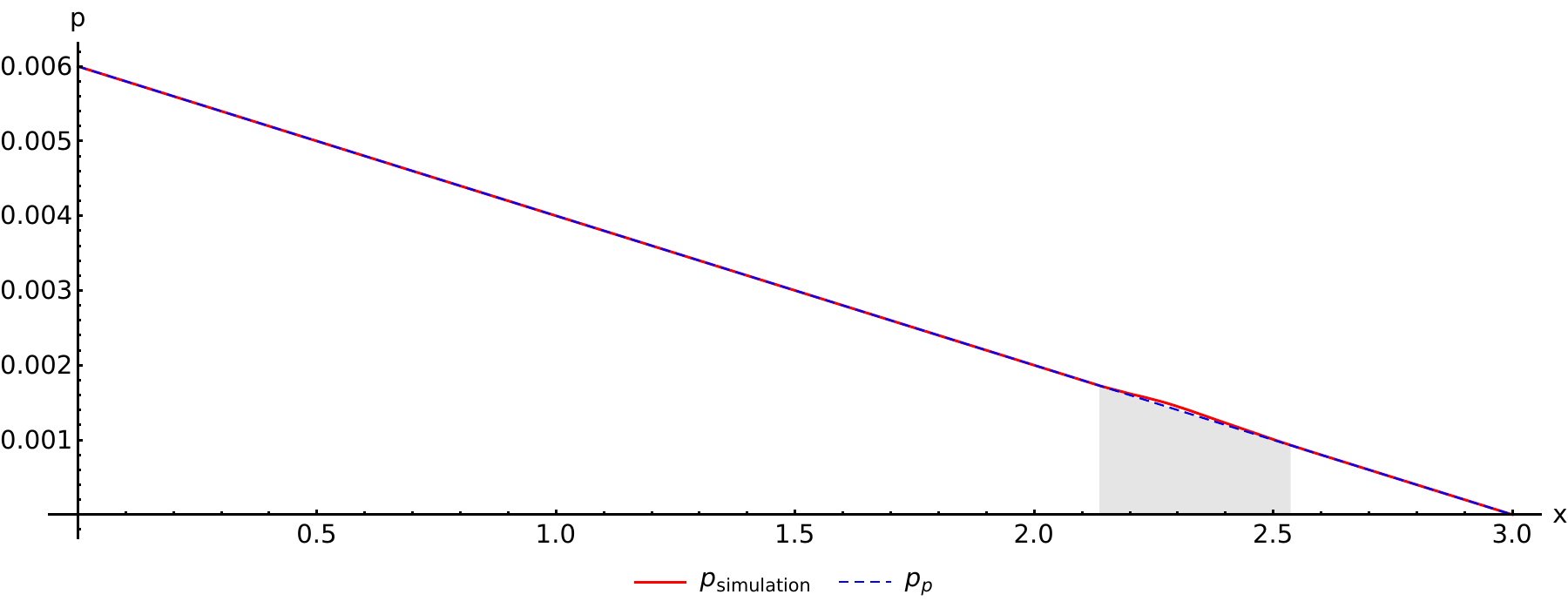}}
    \caption{Cross section plots of velocity in $y$-direction and pressure in $x$-direction and Poiseuille profile at $T = 40$ with gray-shaded rigid body domain.}
    \label{fig:poiseuille_results}
\end{figure}

We compute fully discrete approximations in $(\mathbb P_1(I_n;H_h^2))^2 \times \mathbb P_1(I_n;H_h^1)$, for $n=1,\ldots, N$, for the time step size $\tau = 0.1$ and a background mesh with 223\,750 space-time degrees of freedom in each time step. The mesh and the computed flow profile at the final simulation time $T = 40$ are illustrated in Fig.~\ref{fig:poiseuille_setting}. This solution is fully converged. In fact, the Poiseuille profile is nicely recovered by the simulation. Fig.~\ref{fig:poiseuille_results} shows the cross sections plots in $y$ and $x$ coordinate directions, respectively, of the flow and pressure profile computed by CutFEM. The position of the rigid body domain is gray-shaded. In this domain the velocity and pressure values are computed by the extension and ghost penalty stabilization operator defined in \eqref{eq:extension_operator}. The computed profiles and the Poiseuille profile nicely coincide which clearly demonstrates the accuracy of the suggested CutFEM approach and efficiency of the combined extension and stabilization \eqref{eq:extension_operator}.  

\subsection{Dynamic flow around moving cylinder}
\label{sec:dynamic_flow_around_moving_cylinder}

\begin{figure}[h!t]
	\centering
    \subcaptionbox{Velocity field with overlay of moving rigid domain. \label{fig:moving_ball_v_1}}
	[0.49\columnwidth]{
    \includegraphics[width=0.45\textwidth,keepaspectratio]{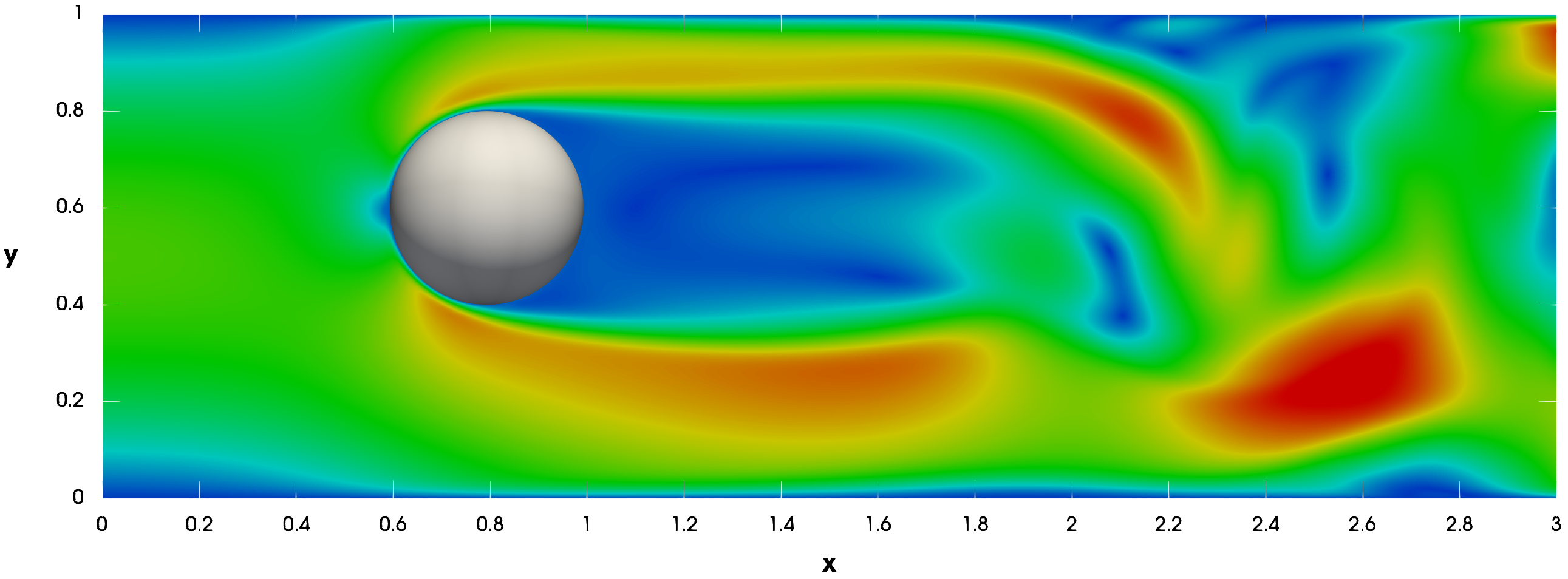}
	}
    \subcaptionbox{Velocity field with computational background mesh. \label{fig:moving_ball_v_grid_1}}
	[0.49\columnwidth]{
	\includegraphics[width=0.45\textwidth,keepaspectratio]{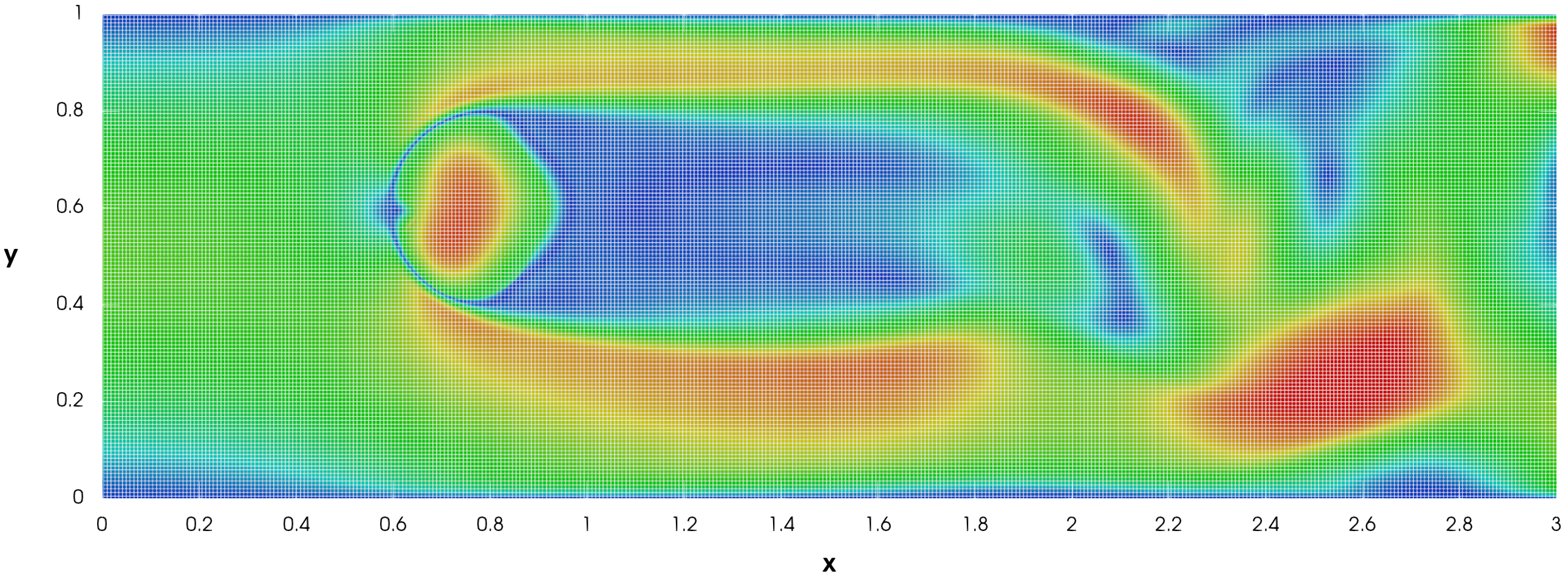}
	}
    \\[2ex]
    \subcaptionbox{Pressure field with overlay of moving rigid body. \label{fig:moving_ball_p_1}}
	[0.49\columnwidth]{
    \includegraphics[width=0.45\textwidth,keepaspectratio]{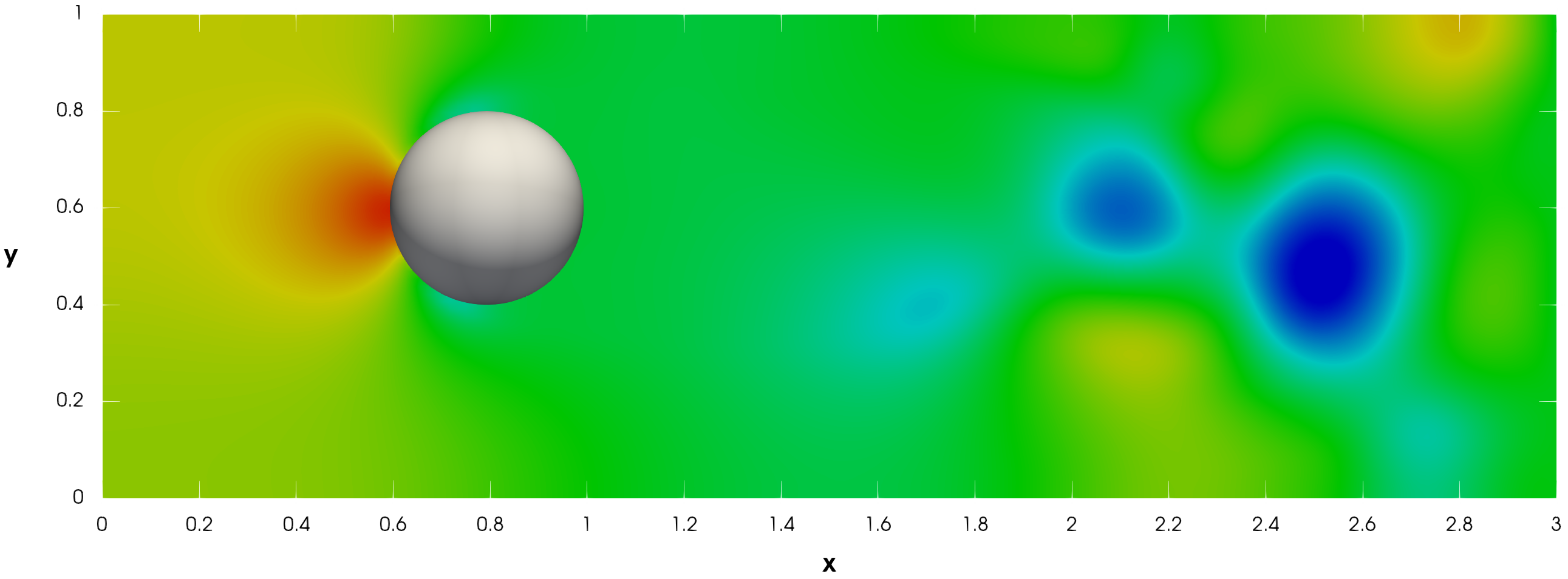}
	}
    \subcaptionbox{Pressure field with computational background mesh. \label{fig:moving_ball_p_grid_1}}
	[0.49\columnwidth]{
	\includegraphics[width=0.45\textwidth,keepaspectratio]{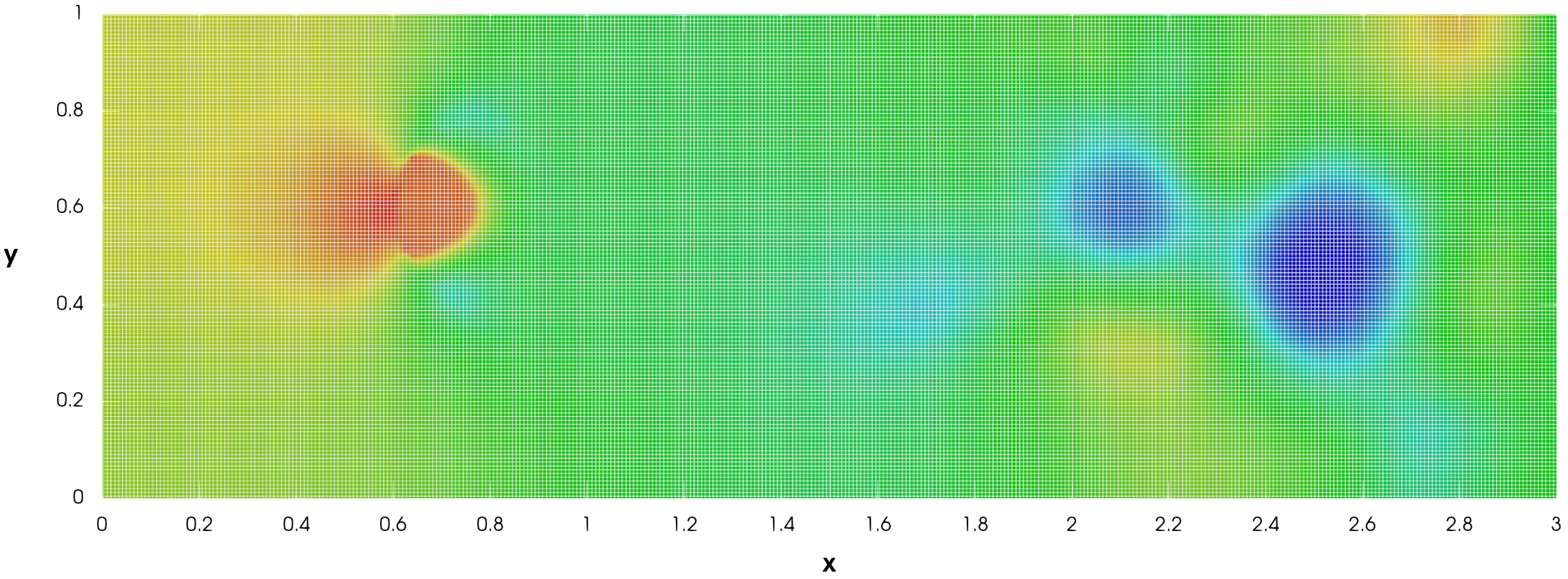}
	}
    \caption{Solution of the CutFEM simulation at time $t = 10.12$.}
    \label{fig:moving_plots_1}
    \vspace*{3ex}
	\centering
    \subcaptionbox{Velocity field with overlay of moving rigid body. \label{fig:moving_ball_v_2}}
	[0.49\columnwidth]{
    \includegraphics[width=0.45\textwidth,keepaspectratio]{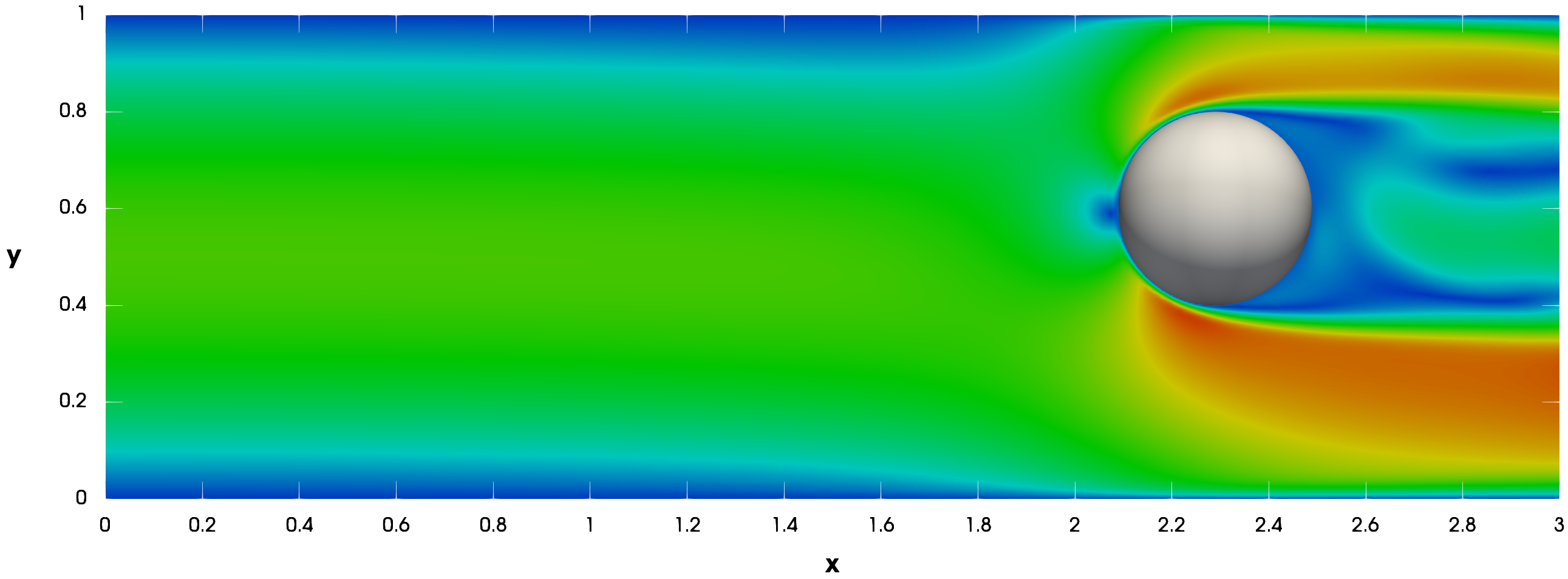}
	}
    \subcaptionbox{Velocity field with computational background mesh. \label{fig:moving_ball_v_grid_2}}
	[0.49\columnwidth]{
	\includegraphics[width=0.45\textwidth,keepaspectratio]{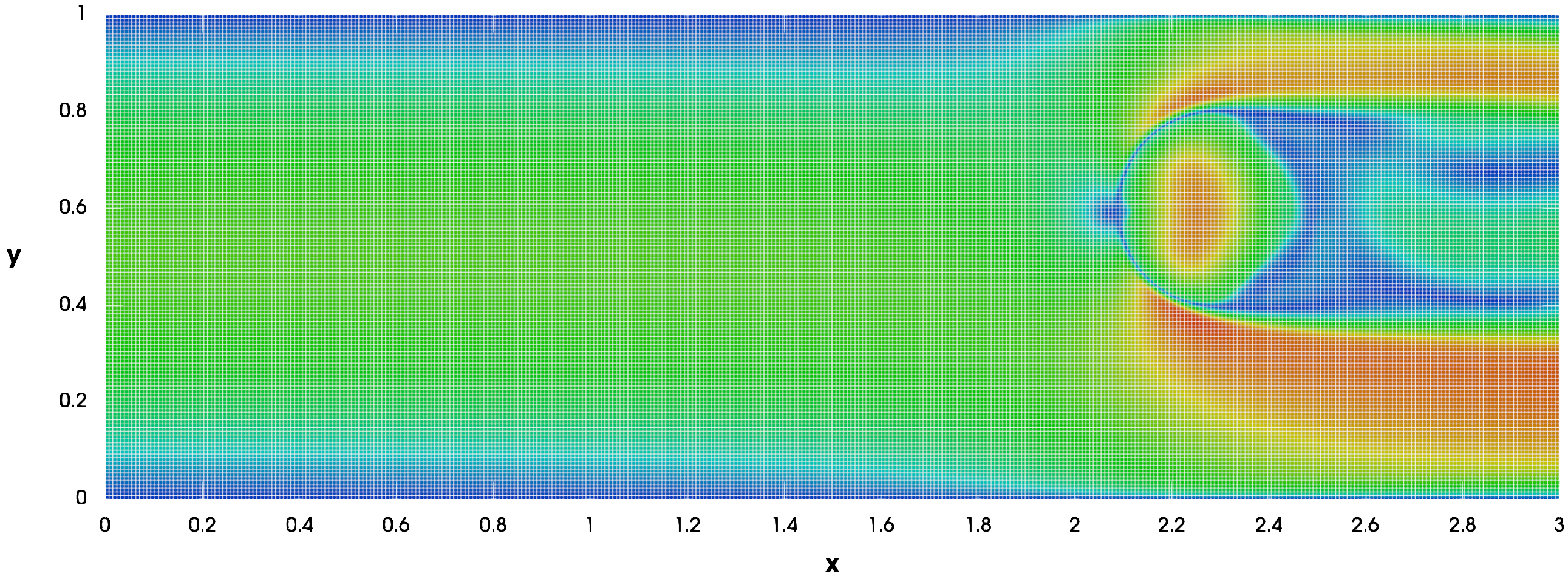}
	}
    \\[2ex]
    \subcaptionbox{Pressure field with overlay of rigid body.
    \label{fig:moving_ball_p_2}}
	[0.49\columnwidth]{
    \includegraphics[width=0.45\textwidth,keepaspectratio]{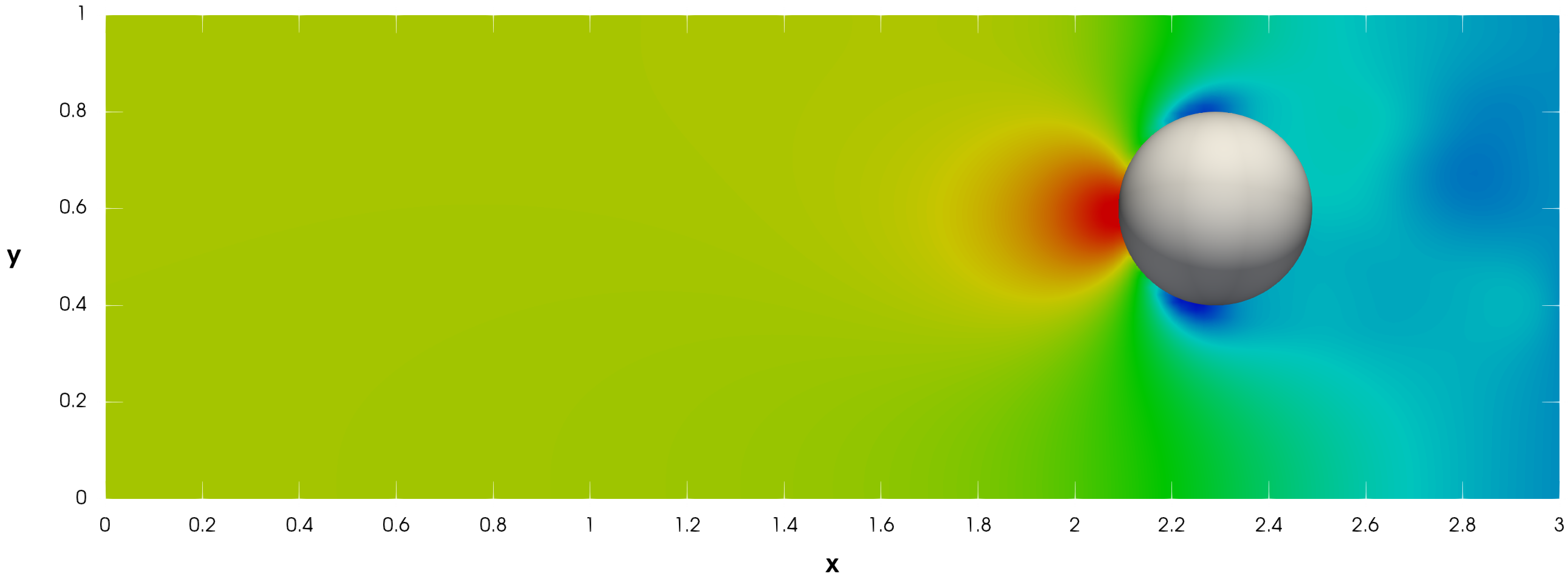}
	}
    \subcaptionbox{Pressure field with computational background mesh. \label{fig:moving_ball_p_grid_2}}
	[0.49\columnwidth]{
	\includegraphics[width=0.45\textwidth,keepaspectratio]{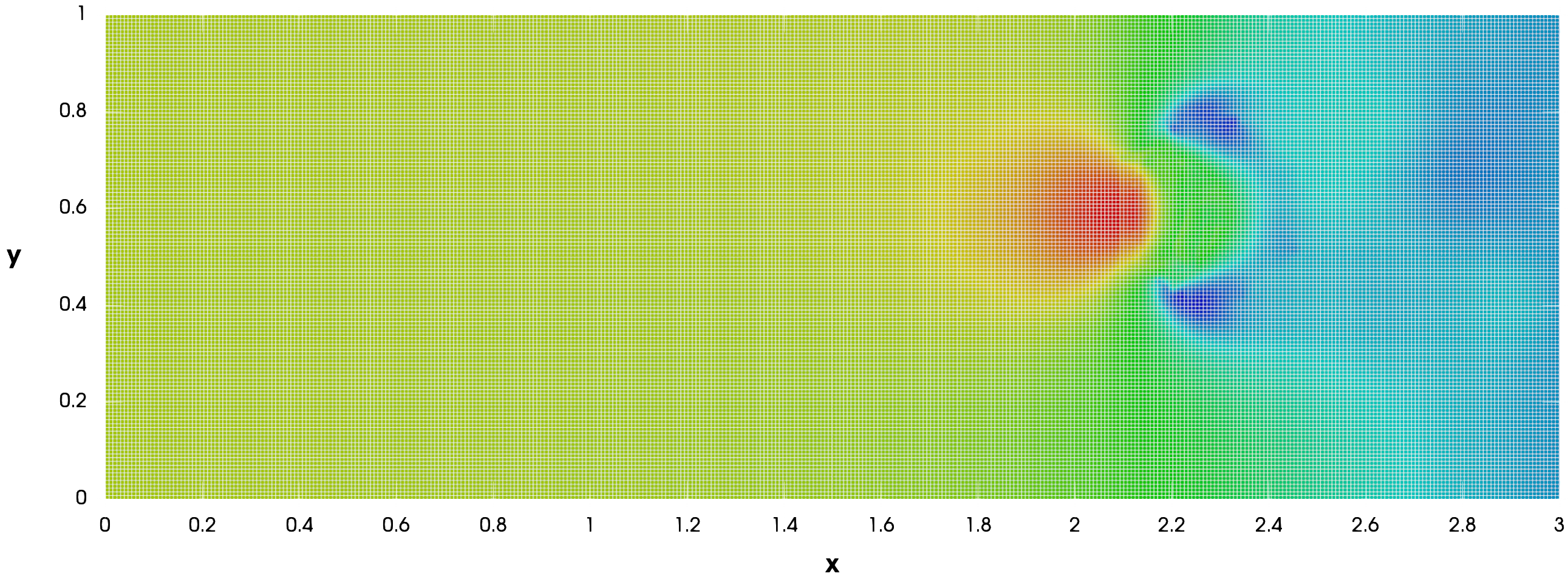}
	}
	\\[2ex]
	\subcaptionbox{Scaling for the velocity and pressure plots. \label{fig:pressure_scale}}
	[0.97\columnwidth]{
	\includegraphics[width=0.4\textwidth,keepaspectratio]{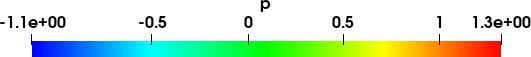}
	\hspace*{8ex}
	\includegraphics[width=0.4\textwidth,keepaspectratio]{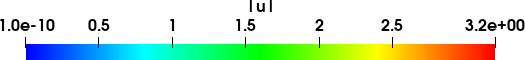}
	}
    \caption{Solution of the CutFEM simulation at $t = T$}
    \label{fig:moving_plots_2}
\end{figure}

In the last numerical example we illustrate the stability and performance properties of the proposed CutFEM for a problem of higher interest in practice. Precisely, we simulate dynamic flow around a moving (two-dimensional) ball. The background domain is $\Omega = (0, 3) \times (0, 1)$, and we let $I = (0, 29]$. At the left inflow boundary we prescribe the parabolic inflow profile
\begin{equation}
\vec{g}_i(\vec{x},t) =
\begin{cases}
\left(6 t y (1 - y),0\right)^\top\,, & t \leq 1|,, \\
\left(6 y (1 - y),0\right)^\top\,, & t > 1\,. \\
\end{cases}
\end{equation}
For the rigid body, a ball with a radius $r = 0.2$ and an initial position of its center at $\vec{x}_r(0) = (1.545, 0.6)^\top$ is chosen. For a non-moving ball, the flow setting would result in a Reynolds number of $Re = 400$. The motion of the ball's center is prescribed by \cref{eq:rigid_motion}, with $A = 0.8$ and $\omega = 0.5$.
On the boundary $\Gamma_r^t$ of the ball, the no-slip condition $\vec v = \vec v_r=(A \omega \cos(\omega t),0)^\top$ for a moving, rigid domain is applied, such that the fluid and ball velocity coincide on their interface, cf.\,\cite[p.\,47]{richterFluidstructureInteractionsModels2017}. We use the time step size $\tau = 0.01$ and a uniform structured background mesh with $h = \frac{1}{64 \sqrt{2}}$, which results in 889\,862 degrees of freedom in each time interval.
With this setup the rigid interface crosses at max one spatial mesh cell.
Discrete solutions are computed in $(\mathbb P_1(I_n;H_h^2))^2 \times \mathbb P_1(I_n;H_h^1)$. Fig.~\ref{fig:moving_plots_1} and Fig.~\ref{fig:moving_plots_2} illustrate the computed profiles of velocity and the pressure at time $t = 10.12$ and $t = T$, respectively.
The figures show that by the application of the ghost penalty stabilization and extension of the discrete solution to the ghost rigid domain, defined by \eqref{eq:extension_operator}, lead to stable approximations of the flow problem on a background mesh. The color legend of Fig.~\ref{fig:moving_plots_2} shows that reasonable values for the velocity and pressure variable are obtained in the ghost subdomain of the rigid ball. We note that all kind of cuts (cf.\ Subsec.~\ref{Subsec:IterInt}), including irregular ones with small portions, arise in the simulation. Unphysical oscillations that are due to irregular cuts of finite elements or insufficient extensions to the rigid domain, are strongly reduced and do not perturb the ambient fluid flow.

\section{Summary and outlook}
\label{Sec:SumOut}

In this work a CutFEM approach for the numerical simulation of fluid flow on evolving domain was presented. The physical domain was embedded into a fixed computational background mesh. Discontinuous Galerkin methods and inf-sup stable pairs of finite element spaces were used for the discretization of the time and space variables, respectively. Ghost penalty stabilization for the treatment of irregular cuts along with an extension of the discrete solution to the ghost rigid domain was suggested further. The numerical performance properties of the parallel implementation of the family of schemes were illustrated. Thereby, the accuracy and stability of the approach was confirmed. Here, a parallel direct linear solver was applied for the sake of simplicity. The integration of cuts cells into our geometric multigrid preconditioner \cite{anselmannGeometricMultigridMethod2021} for GMRES iterations for solving the Newton linearized Navier--Stokes is an ongoing work. In this work, the motion of the rigid body was prescribed explicitly. For the future, more sophisticated applications in that the motion of the body is no longer prescribed explicitly (cf., e.g.,  \cite{freiLocModFELocallyModified2021}) such that the interface has to be computed simultaneously or is given by, for instance fluid-structure interaction, are also in the scope of our interest.

\bibliographystyle{ieeetr}
\bibliography{library} 

\end{document}